\RequirePackage{fix-cm}
\documentclass[onecolumn,natbib]{svjour3}%envcountreset]{svjour3}
\smartqed  % flush right qed marks, e.g. at end of proof

\usepackage{rotating}
\usepackage{graphicx}

\usepackage{empheq}

\usepackage{times}
\usepackage{amstext, amsmath, amssymb}

\usepackage{subfigure}

\usepackage{multirow}
\usepackage{bigstrut}

\usepackage{paralist}

\usepackage{url}
\usepackage[T1]{fontenc}

\usepackage{nameref}

\usepackage[%
    colorlinks=true,%
    bookmarksnumbered=true,%
    bookmarksopen=true,%
    citecolor=blue,%
    urlcolor=blue,%
    unicode=true,           % enable unicode encoded PDF strings
    breaklinks=true,        % allow links to break over lines by making
    pdftitle={State-of-the-Art in Sequential Change-Point Detection},%
    pdfauthor={Aleksey S. Polunchenko and Alexander G. Tartakovsky},%
]{hyperref}

\newcommand{\ignore}[1]{}

\renewcommand{\Pr}{\mathbb{P}} % probability
\DeclareMathOperator{\EV}{\mathbb{E}} % expected value

\DeclareMathOperator{\LR}{\Lambda}

\DeclareMathOperator{\ADD}{ADD}
\DeclareMathOperator{\PFA}{PFA}

\DeclareMathOperator{\RIADD}{RIADD}

\DeclareMathOperator*{\esssup}{ess\,sup}

\newcommand{\tinyinfty}{\infty}

\newcommand{\T}{T}

\renewcommand{\le}{\leqslant} % AMS le ge
\renewcommand{\ge}{\geqslant}

\DeclareMathOperator{\ONE}{\mathchoice{\rm 1\mskip-4.2mu l}{\rm 1\mskip-4.2mu l}{\rm 1\mskip-4.6mu l}{\rm 1\mskip-5.2mu l}}
\newcommand{\indicator}[1]{\ONE_{\left\{#1\right\}}}

\spnewtheorem*{remark*}{Remark}{\bf}{\rm}

\usepackage{color}

\newcommand{\at}[1]{\textcolor{blue}{#1 (AT)}}

\graphicspath{{./gfx/}}

\journalname{Methodol Comput Appl Probab}

\begin{document}

%TITLE PAGE
%The title page should include:
%
%    * The name(s) of the author(s)
%    * A concise and informative title
%    * The affiliation(s) and address(es) of the author(s)
%    * The e-mail address, telephone and fax numbers of the corresponding author
%

\title{State-of-the-Art in Sequential Change-Point Detection\thanks{This work was supported by the U.S.\ Army Research Office under MURI grant  W911NF-06-1-0044, by the U.S.\ Air Force Office of Scientific Research under MURI grant FA9550-10-1-0569, by the U.S.\ Defense Threat Reduction Agency under grant HDTRA1-10-1-0086, and by the U.S.\ National Science Foundation under grants CCF-0830419 and EFRI-1025043.}
}
%\subtitle{The discrete-time case}

\titlerunning{State-of-the-art in change-point detection}  % if too long for running head

\author{Aleksey\ S.\ Polunchenko
    \and
        Alexander\ G.\ Tartakovsky
}

\authorrunning{Polunchenko and Tartakovsky} % if too long for running head

\institute{A.S.\ Polunchenko \at
        Department of Mathematics, University of Southern California,\\
        3620 S. Vermont Ave., KAP 108\\
        Los Angeles, CA 90089-2532, USA\\
        Tel.: +1-213-821-1892 \\
        \email{polunche@usc.edu} % \\
    \and
        A.G.\ Tartakovsky \at
        Department of Mathematics, University of Southern California,\\
        3620 S. Vermont Ave., KAP 108\\
        Los Angeles, CA 90089-2532, USA\\
        Tel.: +1-213-740-2450, Fax: +1-213-740-2424\\
        \email{tartakov@usc.edu}  %  \\
}

\date{Received: date / Accepted: date}

\maketitle % END OF TITLE PAGE

%Abstract
% 150 to 250 words; no undefined abbreviations or unspecified references;
\begin{abstract}
We provide an overview of the state-of-the-art in the area of sequential change-point detection assuming discrete time and known pre- and post-change distributions. The overview spans over all major formulations of the underlying optimization problem, namely, Bayesian, generalized Bayesian, and minimax. We pay particular attention to the latest advances in each. Also, we link together the generalized Bayesian problem with multi-cyclic disorder detection in a stationary regime when the change occurs at a distant time horizon. We conclude with two case studies to illustrate the cutting edge of the field at work.

%Keywords
% provide 4 to 6 keywords which can be used for indexing purposes.

\keywords{
CUSUM chart\and
Quickest change detection\and
Sequential analysis\and
Sequential change-point detection\and
Shiryaev's procedure\and
Shiryaev--Roberts procedure\and
Shiryaev--Roberts--Pollak procedure\and
Shiryaev--Roberts--$r$ procedure
}


%MSC2010
%
%  Statistics -> Sequential methods
%   62L10 - Sequential analysis
%   62L15 - Optimal stopping
%
%  Probability theory and stochastic processes -> Stochastic processes
%   60G40 - Stopping times; optimal stopping problems; gambling theory
%
%  Statistics -> Parametric inference
%   62F12 - Asymptotic properties of estimators
%   62F15 - Bayesian inference
%
%  Statistics -> Decision theory
%   62C10 - Bayesian problems; characterization of Bayes procedures
%   62C20 - Minimax procedures
%
\subclass{MSC 62L10\and MSC 60G40\and MSC 62C10\and MSC 62C20}
\end{abstract}

%-------------------------------------------------------------------------------------------------%
\section{Introduction}
\label{sec:intro} % LABEL NOT REFERENCED, YET TO BE KEPT

Sequential change-point detection (or quickest change detection, or quickest ``disorder'' detection) is concerned with the design and analysis of techniques for {\em quickest} (on-line) detection of a change in the state of a phenomenon, subject to a tolerable limit on the risk of a false detection. Specifically, the substrate of the phenomenon is a time process that may unexpectedly undergo an abrupt change-of-state from ``normal'' to ``abnormal'', each defined as deemed appropriate given the physical context at hand. Inference about the current state of the process is drawn by virtue of (quantitative) observations (e.g., measurements). The {\em sequential setting} assumes the observations are made successively, and, so long as the behavior thereof suggests the process is in the normal state, it is let to continue. However, if  the state is believed to have altered, one's aim is to detect the change ``as soon as possible'', so that an appropriate response can be provided in a timely manner. Thus, with the arrival of every new observation one is faced with the question of whether to let the process continue, or to stop it and raise an alarm (and, e.g., investigate). The decision has to be made in real time based on the available data. The time instance at which the process' state changes is referred to as the {\em change-point}, and the challenge is that it is not known in advance.

Historically, the subject of change-point detection first began to emerge in the 1920--1930's motivated by considerations of quality control. Shewhart's charts were popular in the past (see~\citealp{Shewhart:Book31}). Efficient (optimal and quasi-optimal) sequential detection procedures were developed much later in the 1950-1960's, after the emergence of Sequential Analysis, a branch of statistics ushered by~\cite{Wald:Book47}. The ideas set in motion by Shewhart and Wald have formed a platform for a vast literature on both theory and practice of sequential change-point detection. See, e.g.,~\cite{Girschick+Rubin:AMS52},~\cite{Page:B54},~\cite{Shiryaev:SMD61, Shiryaev:TPA63,Shiryaev:Book78},~\cite{Roberts:T66},~\cite{Siegmund:Book85},~\cite{Tartakovsky:Book91},~\cite{Brodsky+Darkhovsky:Book93},~\cite{Basseville+Nikiforov:Book93},~\cite{Poor+Hadjiliadis:Book08}.

The desire to detect the change quickly causes one to be trigger-happy, which, on one hand, will lead to an unacceptably high level of the risk of sounding a {\em false alarm} -- terminating the process prematurely as a result of an erroneous decision that the change did occur, while, in fact, it never did. On the other hand, attempting to avoid false alarms too strenuously will cause a long delay between the actual time of occurrence of the change (i.e., the true change-point) and the time it is detected. Hence, the essence of the problem is to attain a tradeoff between two contradicting performance measures -- the loss associated with the delay in detection of a true change and that associated with raising a false alarm. A good sequential detection policy is expected to minimize the average loss related to the detection delay, subject to a constraint on the loss associated with false alarms (or vice versa).

\begin{figure}
    \centering
    \includegraphics[width=0.95\textwidth]{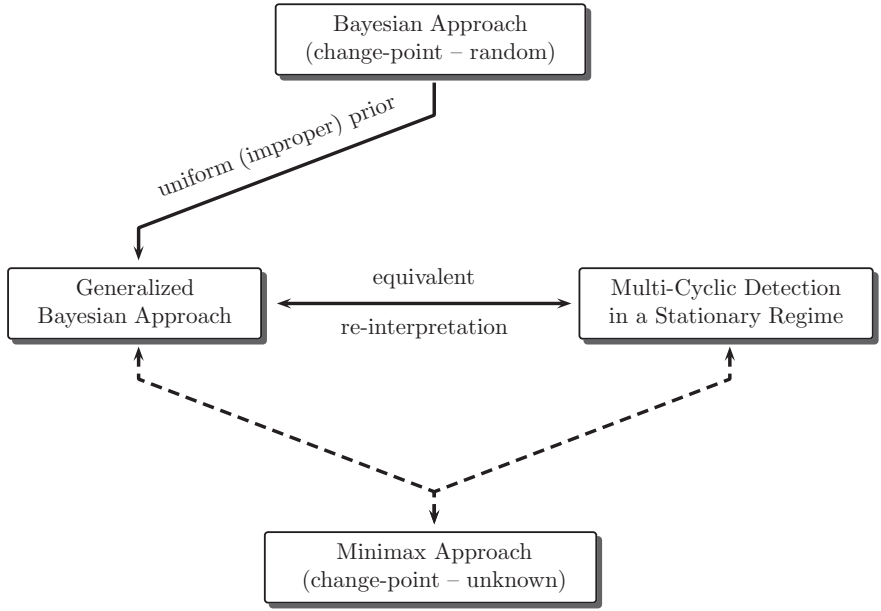}
    % figure caption to be below the figure
    \caption{Four approaches to sequential quickest change-point detection.}
    \label{fig:change-point-detection-big-picture}
\end{figure}

Putting this idea on a rigorous mathematical basis requires formal definition of both the ``detection delay'' and the ``risk of raising a false alarm''. To this end, contemporary theory of sequential change-point detection distinguishes four different approaches: the minimax approach, the Bayesian approach, the generalized Bayesian approach, and the approach related to multi-cyclic detection of a distant change in a stationary regime. Alone, each has its own history and area(s) of application. This notwithstanding the four approaches are connected and fit together into one big picture shown in~\autoref{fig:change-point-detection-big-picture}.

The aim of this paper is to give a brief {\it expos\'{e}} of the above four approaches to quickest change detection. Specifically, the plan is to assess the progress made to date within each with the emphasis on the novel exact and asymptotic optimality results.

%-------------------------------------------------------------------------------------------------%
\section{Change-point models}
\label{sec:change-point-models} % LABEL NOT REFERENCED

To formally state the general quickest change-point detection problem, we first have to introduce a change-point model as well as a model for the observations. To this end, a myriad of scenarios is possible; see, e.g.,~\cite{Fuh:AS03,Fuh:AS04},~\cite{Lai:JRSS95,Lai:IEEE-IT98},~\cite{Shiryaev:SMD61,Shiryaev:TPA63,Shiryaev:Book78,Shiryaev:TPA09,Shiryaev:SA10},~\cite{Tartakovsky:Book91,Tartakovsky:TPA09},~\cite{Tartakovsky+Moustakides:SA10},~\cite{Tartakovsky+Veeravalli:TPA05}. This section is intended to review the major ones.

A change-point model is characterized by the probabilistic structure of the monitored process (independent, identically or non-identically distributed, correlated, etc.) as well as by that of the change-point (unknown deterministic, random completely or partially dependent on the observed data, random fully independent from the observations).

Consider a probability triple $(\Omega,\mathcal{F},\Pr)$, where $\mathcal{F}=\vee_{n\ge 0}\mathcal{F}_n$, $\mathcal{F}_n$ is the sigma-algebra generated by the first $n\ge1$ observations ($\mathcal{F}_{0}=\{\varnothing,\Omega\}$ is the trivial sigma-algebra), and $\Pr\colon\mathcal{F}\mapsto[0,1]$ is a probability measure. Let $\Pr_{\tinyinfty}$ and $\Pr_{0}$ be two mutually locally absolutely continuous (i.e., equivalent) probability measures; for a general case permitting singular measures to be present, see~\cite{Shiryaev:TPA09}. For $d=\{0,\infty\}$, write $\Pr_{d}^{(n)}=\Pr_{d}|_{\mathcal{F}_n}$ for the restriction of measure $\Pr_{d}$ to the  sigma-algebra $\mathcal{F}_n$, and let $p_{d}^{(n)}(\cdot)$ be the density of $\Pr_{d}^{(n)}$ (with respect to a dominating sigma-finite measure).

Let $\{X_n\}_{n\ge1}$ denote the series of (random) observations; the series is defined on the probability space $(\Omega,\mathcal{F},\Pr)$, and distribution-wise is such that for some time index $\nu$, the observations $X_1,X_2,\ldots,X_{\nu}$ adhere to measure $\Pr_{\tinyinfty}$ (``normal'' regime), but $X_{\nu+1},X_{\nu+2},\ldots$ follow measure $\Pr_{0}$ (``abnormal'' regime). That is, at an unknown time instant $\nu$ (change-point), the observations undergo a change-of-regime from normal to abnormal.  Note that $\nu$ is the serial number of the last normal observation, so that if $\nu=0$, then the entire series $\{X_n\}_{n\ge1}$ is in the abnormal regime admitting measure $\Pr_{0}$, while if $\nu=\infty$, then $\{X_n\}_{n\ge1}$ is in the normal regime admitting measure $\Pr_{\infty}$ (i.e., there is no change) .  Another practice popular in the literature is to define $\nu$ as the serial number of the first post-change observation. Although the two definitions map into one another, throughout the remainder of the paper we will follow the former convention.

For every fixed $\nu\ge0$, the change-of-regime in the series $\{X_n\}_{n\ge1}$ gives rise to a new probability measure $\Pr_{\nu}$.  We will now demonstrate how to construct the pdf $p_{\nu}^{(n)}(\boldsymbol{X}_1^n)$ of $\Pr_{\nu}^{(n)}$ for $n\ge1$ and $\nu\ge0$ in the most general case. For the sake of brevity, we will omit the superscript and will write $p_{\nu}(\boldsymbol{X}_1^n)$ in the following.

For $1\le i\le j$, let $\boldsymbol{X}_i^j=(X_i,X_{i+1},\ldots,X_j)$, that is, $\boldsymbol{X}_i^j$ is a sample of $j-i+1$ successive observations indexed from $i$ through $j$. Hence, if the sample $\boldsymbol{X}_1^n=(X_1,X_{2},\ldots,X_n)$ is observed, then $\boldsymbol{X}_1^k=(X_1,\ldots,X_k)$ is the vector of the first $k$ observations in this sample and $\boldsymbol{X}_{k+1}^n=(X_{k+1},\ldots,X_n)$ is the vector of the rest of the observations in the sample from $k+1$ to $n$.

First, suppose $\nu$ is deterministic unknown. This is the main assumption of the minimax approach; recall~\autoref{fig:change-point-detection-big-picture}. To get density $p_{\nu}(\boldsymbol{X}_1^n)$, observe that by the Bayes rule
\begin{align*}
p_{\tinyinfty}(\boldsymbol{X}_1^n)&=
p_{\tinyinfty}(\boldsymbol{X}_1^{\nu})\times
p_{\tinyinfty}(\boldsymbol{X}_{\nu+1}^n|\boldsymbol{X}_1^{\nu})
\;\;\text{and}\;\;
p_{0}(\boldsymbol{X}_1^n)=
p_{0}(\boldsymbol{X}_1^{\nu})\times
p_{0}(\boldsymbol{X}_{\nu+1}^n|\boldsymbol{X}_1^{\nu}),
\end{align*}
whence by combining the first factor of the pre-change density, $p_{\tinyinfty}(\boldsymbol{X}_1^n)$, with the second one of the post-change density, $p_{0}(\boldsymbol{X}_1^n)$, we obtain $p_{\nu}(\boldsymbol{X}_1^n)=p_{\tinyinfty}(\boldsymbol{X}_1^{\nu})\times p_{0}(\boldsymbol{X}_{\nu+1}^n|\boldsymbol{X}_1^{\nu})$,
or, after some more algebra using the Bayes rule,
\begin{align}\label{eq:non-iid-model-def}
p_{\nu}(\boldsymbol{X}_1^n)
&=
\left(\,\prod_{j=1}^{\nu} p_{\tinyinfty}^{(j)}(X_j|\boldsymbol{X}_1^{j-1})\right)
\times
\left(\,\prod_{j=\nu+1}^n p_{0}^{(j)}(X_j|\boldsymbol{X}_1^{j-1})\right),
\end{align}
where $p_{\tinyinfty}^{(j)}(X_j|\boldsymbol{X}_1^{j-1})$ and $p_{0}^{(j)}(X_j|\boldsymbol{X}_1^{j-1})$ are the conditional densities of the $j$-th observation, $X_j$, given the past information $\boldsymbol{X}_1^{j-1}$, $j\ge1$. Note that in general these densities depend on $j$. Hereafter it is understood that $\prod_{j=k+1}^n p_{d}^{(j)}(X_j | \boldsymbol{X}_1^{j-1})=1$ for $k\ge n$.

Model~\eqref{eq:non-iid-model-def} is very general. It does not assume either independence or homogeneity of observations --- the observations may be arbitrary dependent and nonidentically distributed. Furthermore, in certain state-space models and hidden Markov models due to the propagation of the change-point the post-change conditional densities $p_{0}^{(j, \nu)} (X_j|\boldsymbol{X}_1^{j-1})$, $\nu+1\le j\le n$ depend on the change-point $\nu$; see, e.g.,~\cite{Tartakovsky:TPA09}. %and~\cite{Fuh:AS03,Fuh:AS04}.
Model~\eqref{eq:non-iid-model-def} includes practically all possible scenarios. If, for example, there is a switch of one non-iid model to another non-iid model, which are  mutually independent,  then the two segments, pre- and post-change, of the observed process are independent, and in \eqref{eq:non-iid-model-def} the post-change conditional densities $p_{0}^{(j)} (X_j|\boldsymbol{X}_1^{j-1})$, $j \ge \nu+1$ are replaced with $p_{0}^{(j)} (X_j|\boldsymbol{X}_1^{\nu})$.

Suppose now that the observations $\{X_n\}_{n\ge1}$ are {\em independent} and such that $X_1,\ldots,X_{\nu}$ are each distributed according to a common density $f(x)$, while $X_{\nu+1},X_{\nu+2},\ldots$ each follow a common density $g(x)\not\equiv f(x)$. This is the simplest and most prevalent case. For convenience, from now on it will be referred to as the {\em iid case}, or the {\em iid model}. It can be seen that in this case, model~\eqref{eq:non-iid-model-def} reduces to
\begin{align}\label{iidmodel}
p_{\nu}(\boldsymbol{X}_1^n)
&=
\left(\,\prod_{j=1}^{\nu} f(X_j)\right)
\times
\left(\,\prod_{j=\nu+1}^n g(X_j)\right),
\end{align}
and it will be referenced repeatedly throughout the paper.

If the change-point, $\nu$, is random, which is the ground assumption of the Bayesian approach (see~\autoref{fig:change-point-detection-big-picture}), then the model has to be supplied with the change-point's {\em prior distribution}. There may be several change-point mechanisms and, as a result,  a random variable $\nu$ may be dependent on the observations or independent from the observations. To account for these possibilities at once, let $\pi_{0}=\Pr(\nu\le0)$ and $\pi_n=\Pr(\nu=n|\boldsymbol{X}_1^n)$, $n\ge1$, and observe that the series $\{\pi_n\}_{n\ge0}$ is $\{\mathcal{F}_n\}$-adapted. That is, the probability of the change occurring at time instance $\nu=k$ depends on $\boldsymbol{X}_1^k$, the observations' history accumulated up to (and including) time moment $k\ge1$. With the so defined prior distribution one can describe very general change-point models, including those that assume $\nu$ is a $\{\mathcal{F}_n\}$-adapted stopping time; see~\cite{Moustakides:AS08}.

To conclude this section, we note that when the probability series $\{\pi_n\}_{n\ge0}$ depends on the observed data $\{X_n\}_{n\ge1}$, it is argumentative whether $\{\pi_n\}_{n\ge0}$ can be referred to as the change-point's {\em prior} distribution: it can just as well be viewed as the change-point's {\em a posteriori} distribution. However, a deeper discussion of this subject is out of scope to this paper, and from now on, we will assume that $\{\pi_n\}_{n\ge0}$ do not depend on $\{X_n\}_{n\ge1}$, in which case it represents the ``true'' prior distribution.

%-------------------------------------------------------------------------------------------------%
\section{Overview of optimality criteria}
\label{sec:opt-criteria-overview} % NOT REFERENCED

Contemporary theory of change-point detection is an ensemble of the Bayesian approach, the generalized Bayesian approach, the minimax approach, and the approach related to multi-cyclic detection of a disorder in a stationary regime; see~\autoref{fig:change-point-detection-big-picture}. The object of this section is to briefly discuss each problem setting.

A {\em sequential detection procedure} is a stopping time $\T$ adapted to the filtration $\{\mathcal{F}_n\}_{n\ge0}$ induced by the observations $\{X_n\}_{n\ge1}$, i.e., the event $\{\T\le n\}\in\mathcal{F}_n$ for every $n\ge0$. Therefore, after observing $X_1,\ldots,X_{\T}$ it is declared that the change is in effect. That may or may not be the case. If it is not, then $\T\le\nu$, and it is said that a {\em false alarm} has been sounded. Also, note that since $\mathcal{F}_{0}$ is the trivial sigma-algebra, any $\{\mathcal{F}_n\}$-adapted stopping time $\T$ is either strictly positive with probability (w.p.) $1$, or $\T=0$ w.p. $1$. The latter case is clearly degenerate, and to preclude it, from now on we shall assume $\T>0$ w.p. $1$.
\begin{figure}
    \centering
    \subfigure[An example of the behavior of a phenomenon (process) of interest as exhibited through the sequence of observations $\{X_n\}_{n\ge1}$.]{\label{fig:single-run-idea:data}
        \includegraphics[width=0.9\textwidth]{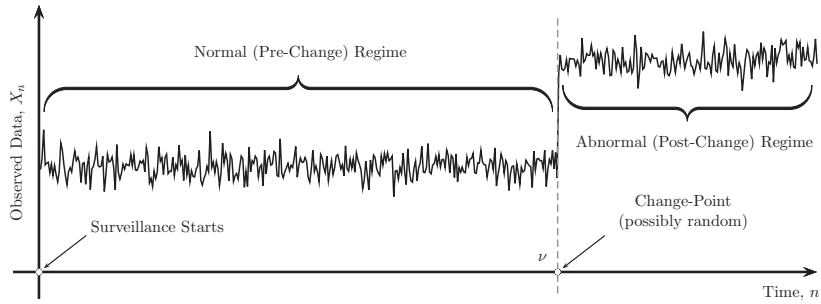}
    } % /subfigure
    \subfigure[Two possible scenarios of the corresponding detection process: false alarm (red trajectory) and correct detection (blue trajectory).]{\label{fig:single-run-idea:detection}
        \includegraphics[width=0.9\textwidth]{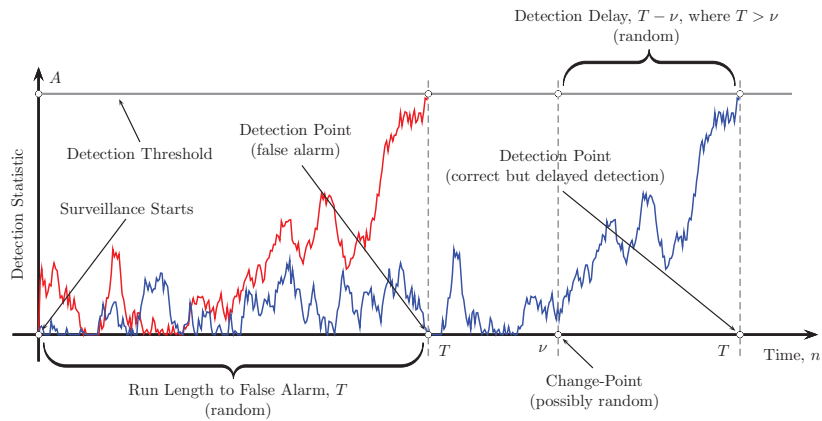}
    } % /subfigure
    \caption{Single-run sequential change-point detection.}
    \label{fig:single-run-idea}
\end{figure}

Common to the Bayesian, generalized Bayesian, and minimax approaches is that the detection procedure is applied {\em only once}; the result is either a false alarm, or a correct (may be delayed) detection. Irrespectively, what takes place beyond the stopping point $\T$ is of no concern. We will refer to this as the single-run paradigm, which is shown in~\autoref{fig:single-run-idea}. Figure~\ref{fig:single-run-idea:data} shows an example of the behavior of a certain process of interest as exhibited through the sequence of observations $\{X_n\}_{n\ge1}$. It can be seen that the process undergoes a shift in the mean at some time instant $\nu$, the change-point. Figure~\ref{fig:single-run-idea:detection} (red trajectory) gives an example of the corresponding detection statistic trajectory that exceeds the detection  threshold prematurely, i.e., before the change occurs. This is a false alarm situation, and $\T$ can be regarded as the (random) run length to the false alarm. Another possibility is shown in Figure~\ref{fig:single-run-idea:detection} (blue trajectory). This is an example where the detection statistic exceeds the detection threshold past the change-point. Note that the detection delay, captured by the difference $\T-\nu$, is random.

In a variety of surveillance applications the detection procedure should be applied {\em repeatedly}. This requires specification of a renewal mechanism after each alarm (false or true). The simplest renewal strategy is to restart from scratch, in which case the procedure becomes multi-cyclic with similar cycles (in a statistical sense) if the process is homogeneous. In the following sections, we will consider such an approach related to detection of a distant change in a stationary regime, assuming that the detection procedure is applied repeatedly starting anew after each time the detection statistic exceeds the threshold.

%-------------------------------------------------------------------------------------------------%
\subsection{Bayesian formulation}

The signature feature of the Bayesian formulation is the assumption that the change-point is a random variable possessing a prior distribution. This is instrumental in certain applications (see~\citealp{Shiryaev:METC06,Shiryaev:SA10}, or~\citealp{Tartakovsky+Veeravalli:TPA05}), but mostly of interest since the limiting versions of Bayesian solutions lead to useful procedures, which are optimal or asymptotically optimal in more practical minimax problems.

Let $\{\pi_k\}_{k\ge0}$ be the prior distribution of the change-point, $\nu$,  where $\pi_{0}=\Pr(\nu\le0)$ and $\pi_k=\Pr(\nu=k)$ for $k \ge 1$.
From the Bayesian point of view, the risk of sounding a false alarm is reasonable to measure by the Probability of False Alarm (PFA), which is defined as
\begin{align}\label{PFAdef}
\PFA^{\pi}(\T)
&=
\Pr^\pi(\T\le\nu)=\sum_{k=1}^\infty \pi_k \Pr_k(T \le k),
\end{align}
where $\Pr^\pi(\mathcal{A}) = \sum\limits_{k=0}^\infty \pi_k \Pr_k(\mathcal{A})$ and the $\pi$ in the superscript emphasizes the dependence on the prior distribution.  Note that summation in \eqref{PFAdef} is over $k \ge 1$ since by convention  $\Pr_k(T \ge 1)=1$, so that $\Pr_k(T \le 0)=0$.
The most popular and practically reasonable way to benchmark the detection delay is through the Average Detection Delay (ADD), which is defined as
\begin{align}\label{ADDBayes}
\ADD^{\pi}(\T)
&=
\EV^\pi [\T-\nu|\T>\nu]=\EV^\pi[(\T-\nu)^+]/\Pr^\pi(\T>\nu),
\end{align}
where hereafter $x^+=\max\{0,x\}$ and $\EV^\pi$ denotes expectation with respect to $\Pr^\pi$.

We are now in a position to formally introduce the notion of Bayesian optimality. Let $\Delta_\alpha=\{\T\colon\PFA^{\pi}(\T)\le\alpha\} $ be the class of detection procedures (stopping times) for which the PFA does not exceed a preset (desired) level $\alpha\in(0,1)$. Then under the Bayesian approach one's aim is to
\begin{align}\label{eq:Bayesian-opt-problem-def}
&\text{find $\T_{\mathrm{opt}}\in\Delta_{\alpha}$ such that $\ADD^{\pi}(\T_{\mathrm{opt}})=
\inf_{\T\in\Delta_{\alpha}}\ADD^{\pi}(\T)$ for every $\alpha\in(0,1)$}.
\end{align}

For the iid model~\eqref{iidmodel} and under the assumption that the change-point $\nu$ has a {\em geometric} prior distribution this problem was solved by~\cite{Shiryaev:SMD61,Shiryaev:TPA63,Shiryaev:Book78}. Specifically, Shiryaev assumed that $\nu$ is distributed according to the zero-modified geometric distribution
\begin{align}\label{eq:Shiryaev-geometric-prior-def}
\Pr(\nu<0)
&=
\pi
\;\;\text{and}\;\;
\Pr(\nu=n)=(1-\pi)p(1-p)^n,
\;\;n\ge0,
\end{align}
where $\pi\in[0,1)$ and $p\in(0,1)$. This is equivalent to choosing the series $\{\pi_n\}_{n\ge0}$ as $\pi_{0}=\Pr(\nu\le0)=\pi+(1-\pi)p$ and $\pi_n=\Pr(\nu=n)=(1-\pi)p(1-p)^n$, $n\ge1$.

Observe now that if $\alpha\ge1-\pi$, then problem~\eqref{eq:Bayesian-opt-problem-def} can be solved by simply stopping right away. This clearly is a trivial solution, since for this strategy the ADD is exactly zero, and $\PFA^{\pi}(\T)=\Pr(\nu>0)=1-\pi$, so that the constraint $\PFA^{\pi}(\T)\le\alpha$ is satisfied.  Therefore, to avoid trivialities we have to assume that $\alpha<1-\pi$. In this case, \cite{Shiryaev:SMD61,Shiryaev:TPA63,Shiryaev:Book78} proved that the optimal detection procedure is based on testing the posterior probability of the change currently being in effect, $\Pr(\nu<n|\mathcal{F}_n)$,  against a certain detection threshold. The procedure stops as soon as $\Pr(\nu<n|\mathcal{F}_n)$ exceed the threshold. This strategy is known as the Shiryaev procedure. To guarantee its strict optimality the detection threshold should be set so as to guarantee that the PFA is exactly equal to the selected level $\alpha$, which is rarely possible.

The Shiryaev procedure will play an important role in the sequel when considering non-Bayes criteria as well. It is more convenient to express Shiryaev's procedure through the average likelihood ratio (LR) statistic
\begin{align}\label{eq:Rnp-S-def}
R_{n,p}
&=
\frac{\pi}{(1-\pi)p}\prod_{j=1}^n\left(\frac{\LR_j}{1-p}\right)+\sum_{k=1}^n\prod_{j=k}^n \left(\frac{\LR_j}{1-p}\right),
\end{align}
where $\LR_n=g(X_n)/f(X_n)$ is the ``instantaneous'' LR for the $n$-th data point, $X_n$. Indeed, by using the Bayes rule, one can show that
\begin{align}\label{eq:posterior-pr-formula}
\Pr(\nu<n|\mathcal{F}_n)
&=
\frac{R_{n,p}}{R_{n,p}+1/p},
\end{align}
whence it is readily seen that ``thresholding'' the posterior probability $\Pr(\nu<n|\mathcal{F}_n)$ is the same as ``thresholding'' the process $\{R_{n,p}\}_{n\ge1}$.  Therefore, the Shiryaev detection procedure has the form
\begin{align}\label{eq:Shirst}
T_{\mathrm{S}}(A)
&=
\inf\{n\ge1\colon R_{n,p}\ge A\},
\end{align}
and if $A=A_\alpha$ can be selected in such a way that the PFA is exactly equal to $\alpha$, i.e., $\PFA^{\pi}(\T_{\rm S}(A_\alpha))=\alpha$, then it is strictly optimal in the class $\Delta(\alpha)$, that is,
\begin{align*}
\inf_{\T\in\Delta(\alpha)}\ADD^\pi(\T)
&=
\ADD^\pi(\T_{\mathrm{S}}(A_\alpha))
\;\;\text{for any}\;\;0<\alpha<1-\pi.
\end{align*}
Note that Shiryaev's statistic $R_{n,p}$ can be rewritten in the recursive form
\begin{align}\label{eq:Rnp-S-recurrent-def}
R_{n,p}
&=
(1+R_{n-1,p})\frac{\LR_n}{1-p},
\;\; n\ge 1,
\;\;\text{with}\;\;
R_{0,p}=\frac{\pi}{(1-\pi)p}.
\end{align}

We also note that~\eqref{eq:Rnp-S-def} and~\eqref{eq:posterior-pr-formula} remain true under the geometric prior distribution~\eqref{eq:Shiryaev-geometric-prior-def} even in the general non-iid case~\eqref{eq:non-iid-model-def}, with $\LR_n= g(X_n|\boldsymbol{X}_1^{n-1})/f(X_n|\boldsymbol{X}_1^{n-1})$. However, in order for the recursion~\eqref{eq:Rnp-S-recurrent-def} to hold in this case, $\{\LR_n\}_{n\ge1}$ should be independent of the change-point.

As $p\to0$, where $p$ is the parameter of the geometric prior~\eqref{eq:Shiryaev-geometric-prior-def}, the Shiryaev detection statistic~\eqref{eq:Rnp-S-recurrent-def} converges to what is known as the {\em Shiryaev--Roberts (SR) detection statistic}. The latter is the basis for the so-called {\em SR procedure}. As we will see, the SR procedure is a ``bridge'' between all four different approaches to change-point detection mentioned above.

For a general asymptotic Bayesian change-point detection theory in discrete time that covers practically arbitrary non-iid models and prior distributions, see~\cite{Tartakovsky+Veeravalli:TPA05}. Specifically, this work addresses the Bayesian approach assuming only that the prior distribution is independent of the observations. The overall conclusion made by the authors is two-fold:
\begin{inparaenum}[\itshape a\upshape)]
    \item the Shiryaev procedure is asymptotically (as $\alpha\to0$) optimal in a very broad class of change-point models and prior distributions, and
    \item depending on the behavior of the prior distribution at the right tail, the SR procedure may or may not be asymptotically optimal.
\end{inparaenum}
Specifically, if the tail is exponential, the SR procedure is not asymptotically optimal, though it is asymptotically optimal if the tail is heavy. When the prior distribution is arbitrary and depends on the observations, we are not aware of any strict or asymptotic optimality results.

%-------------------------------------------------------------------------------------------------%
\subsection{Generalized Bayesian formulation}

The generalized Bayesian approach is the limiting case of the Bayesian formulation, presented in the preceding section. Specifically, in the generalized Bayesian approach the change-point $\nu$ is assumed to be a ``generalized'' random variable with a uniform (improper) prior distribution.

First, return to the Bayesian constrained minimization problem~\eqref{eq:Bayesian-opt-problem-def}. Specifically, consider the iid model~\eqref{iidmodel} and assume that the change-point $\nu$ is distributed according to zero-modified geometric distribution~\eqref{eq:Shiryaev-geometric-prior-def}. Then the Shiryaev  procedure defined in~\eqref{eq:Rnp-S-recurrent-def} and~\eqref{eq:Shirst} is optimal if the threshold $A=A_\alpha$ is chosen so that $\PFA^{\pi}(\T_{\rm S}(A_\alpha))=\alpha$. Suppose now that $\pi=0$ and $p\to0$; this is turning the geometric prior~\eqref{eq:Shiryaev-geometric-prior-def} to an improper uniform distribution. It can be seen that in this case $\{R_{n,p}\}_{n\ge0}$ becomes $\{R_{n,0}\}_{n\ge0}$, where $R_{0,0}=0$ and $R_{n,0}=(1+R_{n-1,0})\LR_n$, $n\ge1$ with $\LR_n=g(X_n)/f(X_n)$. The limit $\{R_{n,0}\}_{n\ge0}$ is known as the SR statistic, and is customarily denoted as $\{R_n\}_{n\ge0}$, i.e., $R_n=R_{n,0}$ for all $n\ge0$; in particular, note that $R_0=0$.

Next, when $\pi=0$ and $p\to0$ it can also be shown that
\begin{align}\label{eq:limits}
\frac{\Pr(\T>\nu)}{p}
&\rightarrow\EV_{\tinyinfty}[\T]
\;\;\text{and}\;\;
\frac{\EV[(\T-\nu)^+]}{p}\rightarrow\sum_{k=0}^\infty\EV_k[(\T-k)^+],
\end{align}
where $\T$ is an arbitrary stopping time. As a result, one may conjecture that the SR procedure minimizes the {\em Relative Integral Average Detection Delay} (RIADD)
\begin{align}\label{eq:RIADD-def}
\RIADD(\T)
&=
\frac{\sum_{k=0}^\infty\EV_k[(\T-k)^+]}{\EV_{\tinyinfty}[\T]}
\end{align}
over all detection procedures for which the {\em Average Run Length (ARL) to false alarm}, $\EV_{\tinyinfty}[\T]$, is no less than $\gamma>1$, an {\it a~priori} set level.

Let
\begin{align}\label{eq:class-Delta-ARL-def}
\Delta(\gamma)
&=
\bigl\{\T\colon\EV_{\tinyinfty}[\T]\ge\gamma\bigr\},
\end{align}
be the class of detection procedures (stopping times) for which the ARL to false alarm $\EV_{\tinyinfty}[\T]$ is ``no worse'' than $\gamma>1$. Then under the generalized Bayesian formulation one's goal is to
\begin{align}\label{eq:genBayesproblem}
&\text{find $\T_{\mathrm{opt}}\in\Delta(\gamma)$ such that $\RIADD(\T_{\mathrm{opt}})=
\inf_{\T\in\Delta(\gamma)}\RIADD(\T)$ for every $\gamma>1$}.
\end{align}

We have already hinted that this problem is solved by the SR procedure. This was formally demonstrated by~\cite{Pollak+Tartakovsky:SS09} in the discrete-time iid case, and by~\cite{Shiryaev:TPA63} and~\cite{Feiberg+Shiryaev:SD06} in continuous time for detecting a shift in the mean of a Brownian motion.

We conclude this subsection with two remarks. First, observe that if the assumption $\pi=0$ is replaced with $\pi=rp$, where $r\ge0$ is a fixed number, then, as $p\to0$, the Shiryaev statistic $\{R_{n,p}\}_{n\ge0}$ converges to $\{R_n^r\}_{n\ge0}$, where $R_n^r=(1+R_{n-1}^r)\LR_n$, $n\ge1$ with $R_0^r=r\ge0$. This is the so-called {\em Shiryaev--Roberts--$r$ (SR--$r$) detection statistic}, and it is the basis for the SR--$r$ detection procedure that starts from an arbitrary deterministic point $r$. This procedure is due to~\cite{Moustakidesetal-SS11}. The SR--$r$ procedure possesses certain minimax properties (cf.~\cite{Polunchenko+Tartakovsky:AS10} and~\cite{Tartakovsky+Polunchenko:IWAP10}). We will discuss this procedure at greater length later.

Secondly, though the generalized Bayesian formulation is the limiting (as $p\to0$) case of the Bayesian approach, it may also be equivalently re-interpreted as a completely different approach -- {\em multi-cyclic disorder detection in a stationary regime}. We will address this approach in~\autoref{sec:opt-criteria-overview:multi-cyclic}.

%-------------------------------------------------------------------------------------------------%
\subsection{Minimax formulation}
\label{sec:opt-criteria-overview:minimax}

Contrary to the Bayesian formulation the minimax approach posits that the change-point is an unknown not necessarily random number. Even if it is random its distribution is unknown. The minimax approach has multiple optimality criteria.

First minimax theory is due to~\cite{Lorden:AMS71} who proposed to measure the risk of raising a false alarm by the ARL to false alarm $\EV_{\tinyinfty}[\T]$ (recall the false alarm scenario from Figure~\ref{fig:single-run-idea:detection}). As far as the risk associated with detection delay is concerned, Lorden suggested to use the ``worst-worst-case'' ADD defined as
\begin{align}\label{eq:SADD-Lorden-def}
\mathcal{J}_{\mathrm{L}}(\T)
&=
\sup_{0 \le \nu<\infty}\biggl\{\esssup\EV_{\nu}[(\T-\nu)^+|\mathcal{F}_{\nu}]\biggr\}.
\end{align}
 Lorden's minimax optimization problem seeks to
\begin{align}\label{eq:Lorden-minimax-problem}
&\text{find $\T_{\mathrm{opt}}\in\Delta(\gamma)$ such that $\mathcal{J}_{\mathrm{L}}(\T_{\mathrm{opt}})=
\inf_{\T\in\Delta(\gamma)}\mathcal{J}_{\mathrm{L}}(\T)$ for every $\gamma>1$} ,
\end{align}
where $\Delta(\gamma)$ is the class of detection procedures with the lower bound $\gamma$ on the ARL to false alarm defined in \eqref{eq:class-Delta-ARL-def}.

For the iid scenario~\eqref{iidmodel},~\cite{Lorden:AMS71} showed that Page's~\citeyearpar{Page:B54} Cumulative Sum (CUSUM) procedure is first-order asymptotically minimax as $\gamma\to\infty$. For any $\gamma>1$, this problem was solved by~\cite{Moustakides:AS86}, who showed that CUSUM  is exactly optimal (see also~\cite{Ritov:AS90} who reestablished Moustakides'~\citeyearpar{Moustakides:AS86} finding using a different decision-theoretic argument).

Though the strict $\mathcal{J}_{\mathrm{L}}(\T)$-optimality of the CUSUM procedure is a strong result, it is more natural to construct a procedure that minimizes the average (conditional) detection delay, $\EV_{\nu}[\T-\nu|\T>\nu]$, for all $\nu\ge0$ simultaneously. As no such uniformly optimal procedure is possible,~\cite{Pollak:AS85} suggested to revise Lorden's version of minimax optimality by replacing $\mathcal{J}_{\mathrm{L}}(\T)$ with
\begin{align}\label{eq:SADD-Pollak-def}
\mathcal{J}_{\mathrm{P}}(\T)
&=
\sup_{0\le\nu<\infty}\EV_{\nu}[\T-\nu|\T>\nu],
\end{align}
the worst conditional expected detection delay. Thus, Pollak's version of the minimax optimization problem seeks to
\begin{align}\label{eq:Pollak-minimax-problem}
&\text{find $\T_{\mathrm{opt}}\in\Delta(\gamma)$ such that $\mathcal{J}_{\mathrm{P}}(\T_{\mathrm{opt}})=
\inf_{\T\in\Delta(\gamma)}\mathcal{J}_{\mathrm{P}}(\T)$ for every $\gamma>1$}.
\end{align}

It is our opinion that $\mathcal{J}_{\mathrm{P}}(\T)$ is better suited for practical purposes for two reasons. First,  Lorden's criterion is effectively a double-minimax approach, and therefore, is overly pessimistic in the sense that $\mathcal{J}_{\mathrm{P}}(\T)\le\mathcal{J}_{\mathrm{L}}(\T)$. Second, it is directly connected to the conventional decision theoretic approach --- the optimization problem \eqref{eq:Pollak-minimax-problem} can be solved by finding the least favorable prior distribution. More specifically, since by the general decision theory, the minimax solution corresponds to the (generalized) Bayesian solution with the least favorable prior distribution, it can be shown that $\sup_{\pi}\ADD^{\pi}(\T)=\mathcal{J}_{\mathrm{P}}(\T)$, where $\ADD^\pi(\T)$ is defined in \eqref{ADDBayes}. In addition, unlike Lorden's minimax problem~\eqref{eq:Lorden-minimax-problem}, Pollak's minimax problem~\eqref{eq:Pollak-minimax-problem} is still not solved.  For these reasons, from now on, when considering the minimax approach, we focus on Pollak's supremum ADD measure $\mathcal{J}_{\mathrm{P}}(\T)$. Some light as to the possible solution (in the iid case) is shed in the work of~\cite{Polunchenko+Tartakovsky:AS10},~\cite{Tartakovsky+Polunchenko:IWAP10}, and~\cite{Moustakidesetal-SS11}. A synopsis of the results is given in Sections~\ref{s:SR-r-opt} and~\ref{sec:case-studies}.

We conclude this section with presenting yet another way to gauge the risk of raising a false alarm, namely, by means of the worst local (conditional) probability of sounding a false alarm within a time ``window'' of a given length. As argued by~\cite{Tartakovsky:IEEE-CDC05,Tartakovsky:SA08-discussion}, in many surveillance applications (e.g., target detection) this probability is a better  option than the ARL to false alarm, which is more global. Specifically, let
\begin{align}\label{eq:class-Delta-sup-PFA-window-def}
\Delta_{\alpha}^m
&=
\biggl\{\T\colon\sup_{k\ge0}\Pr_{\tinyinfty}(k<\T\le k+m|\T>k)\le\alpha\biggr\},
\end{align}
be the class of detection procedures for which $\Pr_{\tinyinfty}(k<\T\le k+m|\T>k)$, the conditional probability of raising a false alarm inside a sliding window of $m\ge1$ observations  is ``no worse'' than a certain {\it a~priori} chosen level $\alpha\in(0,1)$. The size of the window $m$ may either be fixed or go to infinity when $\alpha\to0$.

Let $\T$ be the stopping time associated with a generic detection procedure. The appropriateness of the ARL to false alarm $\EV_{\tinyinfty}[\T]$ as an exhaustive measure of the risk of raising a false alarm is questionable, unless the $\Pr_{\tinyinfty}$-distribution of $\T$ is geometric, at least approximately; see~\cite{Tartakovsky:IEEE-CDC05,Tartakovsky:SA08-discussion}. The geometric distribution is characterized entirely by a single parameter, which a) uniquely determines $\EV_{\tinyinfty}[\T]$, and b) is uniquely determined by $\EV_{\tinyinfty}[\T]$. As a result, if $\T$ is geometric, one can evaluate $\Pr_{\tinyinfty}(k<\T\le k+m|\T>k)$ for  any $k\ge0$ (in fact, for all $k\ge0$ at once).

For the iid model~\eqref{iidmodel},~\cite{Pollak+Tartakovsky:TPA09} showed
that under mild assumptions the $\Pr_{\tinyinfty}$-distribution of the stopping times associated with detection schemes from a certain class is asymptotically (as $\gamma\to\infty$) exponential  with parameter $1/\EV_{\tinyinfty}[\T]$; the convergence is in the $L^p$ sense, where $p\ge1$. See also~\cite{Tartakovsky+Pollak+Polunchenko:IWAP08}. The class includes all of the most popular procedures. Hence, for the iid model~\eqref{iidmodel}, the ARL to false alarm is an acceptable measure of the false alarm rate. However, for a general non-iid model this is not necessarily true, which suggests that alternative measures of the false alarm rate are in order.

As argued by~\cite{Tartakovsky:IEEE-CDC05}, in general, $\sup_{k}\Pr_{\tinyinfty}(k<\T\le k+m|\T>k)\le\alpha$ is a {\em stronger} condition than $\EV_{\tinyinfty}[\T]\ge\gamma$. Hence, in general, $\Delta_{\alpha}^m\subset\Delta(\gamma)$. See also~\cite{Tartakovsky:SA09-discussion}. In~\autoref{sec:case-studies} we take the work of~\cite{Polunchenko+Tartakovsky:AS10} and~\cite{Tartakovsky+Polunchenko:IWAP10} one step further and present a procedure that solves the optimization problem~\eqref{eq:Pollak-minimax-problem} in the class~\eqref{eq:class-Delta-sup-PFA-window-def} in a specific example.

%-------------------------------------------------------------------------------------------------%
\subsection{Multi-cyclic detection of a disorder in a stationary regime}
\label{sec:opt-criteria-overview:multi-cyclic}

Common to all of the above approaches is that the detection procedure is applied only once. This is the single-run paradigm. The result is either a correct (though usually delayed) detection or a false alert; recall~\autoref{fig:single-run-idea}. Yet another formulation may be derived by abandoning the single-run paradigm for the {\em multi-run} or the {\em multi-cyclic} one.

\begin{sidewaysfigure}
    \centering
    \subfigure[An example of the behavior of a process of interest as exhibited through the series of observations $\{X_n\}_{n\ge1}$.]{
        \includegraphics[width=0.9\textwidth]{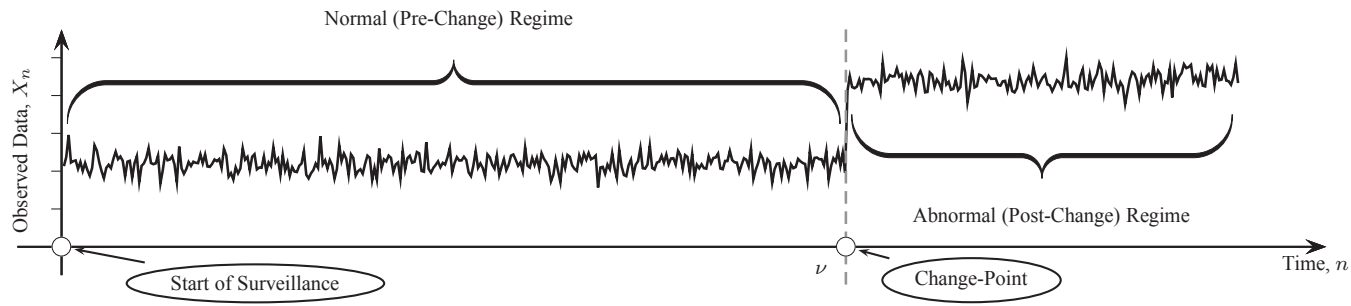}
    } % /subfigure
    \subfigure[An example of the behavior of the detection statistic when the decision to terminate surveillance is made {\em past} the change-point.]{
        \includegraphics[width=0.9\textwidth]{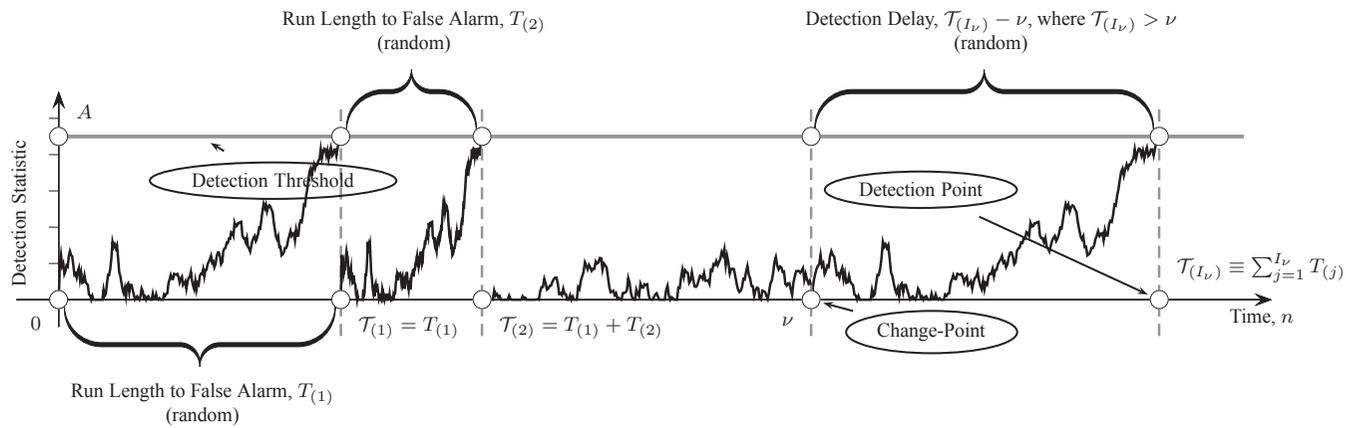}
    } % /subfigure
    % figure caption to be below the figure
    \caption{Multi-cyclic change-point detection in a stationary regime.}
    \label{fig:multi-cyclic-idea}
\end{sidewaysfigure}

Specifically, consider a context in which it is of utmost importance to detect the change as quickly as possible, even at the expense of raising many false alarms (using a repeated application of the same stopping rule) before the change occurs. This is equivalent to saying that the change-point $\nu$ is substantially larger than the tolerable level of false alarms $\gamma$. That is, the change ``strikes'' in a distant future and is preceded by a {\em stationary flow of false alarms}. This scenario is schematically shown in~\autoref{fig:multi-cyclic-idea}. As one can see, the ARL to false alarm in this case is the mean time between (consecutive) false alarms, and therefore may be thought of the false alarm rate (or frequency).

As argued by~\cite{Pollak+Tartakovsky:SS09}, the multi-cyclic approach is instrumental in many surveillance applications, in particular in the areas concerned with intrusion/anomaly detection, e.g., cybersecurity and particularly detection of attacks in computer networks.

Formally, let $T_1,T_2,\ldots$ denote sequential independent repetitions of the same stopping time $\T$, and let $\mathcal{T}_{(j)}=T_{(1)}+T_{(2)}+\cdots+T_{(j)}$ be the time of the $j$-th alarm. Define $I_\nu=\min\{j\ge1\colon\mathcal{T}_{(j)}>\nu\}$. Put otherwise, $\mathcal{T}_{(I_\nu)}$ is the time of detection of the true change that occurs at the time instant $\nu$ after $I_\nu-1$ false alarms have been raised. Write
\begin{align}\label{eq:STADD-def}
\mathcal{J}_{\mathrm{ST}}(\T)
&=
\lim_{\nu\to\infty}\EV_\nu[\mathcal{T}_{(I_\nu)}-\nu]
\end{align}
for the limiting value of the ADD that we will refer to as the {\em stationary ADD} (STADD).

We are now in a position to formalize the notion of optimality in the multi-cyclic setup:
\begin{align}\label{eq:Multicycleoptim}
&\text{find $\T_{\mathrm{opt}}\in\Delta(\gamma)$ such that $\mathcal{J}_{\mathrm{ST}}(\T_{\mathrm{opt}})=
\inf_{\T\in\Delta(\gamma)}\mathcal{J}_{\mathrm{ST}}(\T)$ for every $\gamma>1$}
\end{align}
(among all multi-cyclic procedures).

For the iid model~\eqref{iidmodel}, this problem was solved by~\cite{Pollak+Tartakovsky:SS09}, who showed that the solution is the multi-cyclic SR procedure by arguing that $\mathcal{J}_{\mathrm{ST}}(\T)$ is the same as $\RIADD(\T)$ defined in~\eqref{eq:RIADD-def}. This suggests that the optimal solution of the problem of multi-cyclic change-point detection in a stationary regime is completely equivalent to the solution of the generalized Bayesian problem. The exact result is stated in the next section.

%-------------------------------------------------------------------------------------------------%
\section{Optimality properties of the Shiryaev--Roberts detection procedure}
\label{s:SRoptimality}

From now on we will confine ourselves to the iid scenario~\eqref{iidmodel}, i.e., we assume that a) the observations $\{X_n\}_{n\ge1}$ are independent throughout their history, and b) $X_1,\ldots,X_{\nu}$ are distributed according to a common known pdf $f(x)$ and $X_{\nu+1},X_{\nu+2},\ldots$ are distributed according to a common pdf $g(x)\not\equiv f(x)$, also known.

Let $\mathcal{H}_k\colon\nu=k$ for $\le k<\infty$ and $\mathcal{H}_{\tinyinfty}\colon\nu=\infty$
be, respectively, the hypotheses that the change takes place at the time moment $\nu=k$, $k\ge0$, and that no change ever occurs. The densities of the sample $\boldsymbol{X}_1^n=(X_1,\ldots,X_n)$, $n\ge1$ under these hypotheses are given by
\begin{align*}
\begin{aligned}
p(\boldsymbol{X}_1^n|\mathcal{H}_{\tinyinfty})
&=
\prod_{j=1}^n f(X_j),
\\
p(\boldsymbol{X}_1^n|\mathcal{H}_k)
&=
\prod_{j=1}^k f(X_j)\prod_{j=k+1}^n
g(X_j)\quad\text{for $k<n$},
\end{aligned}
\end{align*}
and $p(\boldsymbol{X}_1^n|\mathcal{H}_{\tinyinfty})=p(\boldsymbol{X}_1^n|\mathcal{H}_k)$ for $k\ge n$, so that the corresponding LR is
\begin{align*}
\LR_n^k
&=
\frac{p(\boldsymbol{X}_1^n|\mathcal{H}_k)}{p(\boldsymbol{X}_1^n|\mathcal{H}_{\tinyinfty})}
=\prod_{j=k+1}^n\LR_j\quad\text{for}\quad k<n,
\end{align*}
where $\LR_n= g(X_n)/f(X_n)$ is the ``instantaneous'' LR for the $n$-th observation $X_n$.

To decide in favor of one of the hypotheses $\mathcal{H}_k$ or $\mathcal{H}_{\tinyinfty}$, the likelihood ratios are then ``fed'' to an appropriate sequential detection procedure, which is chosen according to the particular version of the optimization problem. In this section we are interested in the generalized Bayesian problem stated in~\eqref{eq:genBayesproblem} and in the multi-cyclic disorder detection in a stationary regime stated in~\eqref{eq:Multicycleoptim}. We have already remarked that for the iid model in question the SR procedure solves both these problems. We preface the presentation of the exact results with the introduction of the SR procedure.

%-------------------------------------------------------------------------------------------------%
\subsection{The Shiryaev--Roberts procedure}

The SR procedure is due to the independent work of~\cite{Shiryaev:SMD61,Shiryaev:TPA63} and~\cite{Roberts:T66}. Specifically, Shiryaev considered the problem of detecting a change in the drift of a Brownian motion; Roberts focused on the case of detecting a shift in the mean of an iid Gaussian sequence. The name Shiryaev--Roberts was given by~\cite{Pollak:AS85}, and it has become the convention.

Formally, the SR procedure is defined as the stopping time
\begin{align}\label{eq:T-SR-def}
\mathcal{S}_A
&=
\inf\{n\ge1\colon R_n\ge A\},
\end{align}
where $A>0$ is the detection threshold, and
\begin{align}\label{eq:Rn-SR-def}
R_n
&=
(1+R_{n-1})\LR_n,
\;\;n\ge1\;\;\text{with}\;\; R_0=0
\end{align}
is the SR detection statistic. As usual, we set  $\inf\{\varnothing\}=\infty$, i.e., $\mathcal{S}_A=\infty$ if $R_n$ never crosses $A$.

%-------------------------------------------------------------------------------------------------%
\subsection{Optimality properties}

Recall first that $R_n=\lim_{p\to0}R_{n,p}$, where $R_{n,p}$ is the Shiryaev statistic given by recursion~\eqref{eq:Rnp-S-recurrent-def}. Recall also that the limiting relations~\eqref{eq:limits} hold. These facts allow us to conjecture that the SR procedure is optimal in the generalized Bayesian sense. In addition, as we stated in~\autoref{sec:opt-criteria-overview:multi-cyclic}, the RIADD is equal to the STADD of the multi-cyclic procedure, so that we expect that the repeated SR procedure is optimal for detecting distant changes. The exact result is due to~\cite{Pollak+Tartakovsky:SS09} and is given next.

\begin{theorem}[\citealp{Pollak+Tartakovsky:SS09}]\label{Th:SRoptimality}
Let $\mathcal{S}_A$ be the SR procedure defined by \eqref{eq:T-SR-def} and~\eqref{eq:Rn-SR-def}. Suppose the detection threshold $A=A_{\gamma}$ is selected from the equation $\EV_{\tinyinfty}[\mathcal{S}_{A_{\gamma}}]=\gamma$, where $\gamma>1$ is the desired level of the ARL to false alarm.
\begin{enumerate}[\rm (i)]
\item Then the SR procedure $\mathcal{S}_{A_{\gamma}}$ minimizes $\RIADD(\T)=\sum_{k=0}^\infty\EV_k[(\T-k)^+]/\EV_{\tinyinfty}[\T]$ over all stopping times $\T$ that satisfy $\EV_{\tinyinfty}[\T]\ge\gamma$, that is,
\begin{align*}
\RIADD(\mathcal{S}_{A_{\gamma}})
&=
\inf_{\T\in\Delta(\gamma)}\RIADD(\T)
\;\;\text{for every}\;\;\gamma>1 .
\end{align*}
    \item For any stopping time $\T$, $\RIADD(\T)=\mathcal{J}_{\mathrm{ST}}(\T)$. Therefore, the SR procedure $\mathcal{S}_{A_{\gamma}}$ minimizes the stationary average detection delay among all multi-cyclic procedures in the class $\Delta(\gamma)$, i.e.,
\begin{align*}
\mathcal{J}_{\mathrm{ST}}(\mathcal{S}_{A_{\gamma}})
&=
\inf_{\T\in\Delta(\gamma)}\mathcal{J}_{\mathrm{ST}}(\T)
\;\;\text{for every}\;\;\gamma>1 .
\end{align*}
\end{enumerate}
\end{theorem}

It is worth noting that the ARL to false alarm of the SR procedure satisfies the inequality $\EV_\infty[\mathcal{S}_{A}]\ge A$ for all $A>0$, which can be easily obtained by noticing that $R_n-n$ is a $\Pr_\infty$-martingale with mean zero. Also, asymptotically (as $A\to\infty$),  $\EV_\infty[\mathcal{S}_{A}]\approx  A/\zeta$, where the constant $0< \zeta <1$ is given by~\eqref{kappazeta} below (see~\citealp{Pollak:AS87}).  Hence, setting $A_\gamma=\gamma\zeta$ yields $\EV_\infty[\mathcal{S}_{A_\gamma}]\approx\gamma$, as $\gamma\to\infty$.

%-------------------------------------------------------------------------------------------------%
\section{Optimal and nearly optimal minimax detection procedures}
\label{s:minimax-procedures} % NOT REFERENCED

In this section, we will be concerned exclusively with the minimax problem in Pollak's setting~\eqref{eq:Pollak-minimax-problem}, assuming that the change-point $\nu$ is deterministic unknown. As of today, this problem is not solved in general. As has been indicated earlier, the usual way around this is to consider it asymptotically by allowing the ARL to false alarm $\gamma\to\infty$. The hope is to design such procedure $\T^*\in\Delta(\gamma)$ that $\mathcal{J}_{\mathrm{P}}(\T^*)$ and the (unknown) optimum $\inf_{\T\in\Delta(\gamma)}\mathcal{J}_{\mathrm{P}}(\T)$ will be in some sense ``close'' to each other in the limit, as $\gamma\to\infty$. To this end, the following three different types of asymptotic optimality are usually distinguished.

\begin{definition}[First-Order Asymptotic Optimality]
A procedure $\T^*\in\Delta(\gamma)$ is said to be {\em first-order asymptotically} optimal in the class $\Delta(\gamma)$ if
\begin{align*}
\frac{\mathcal{J}_{\mathrm{P}}(\T^*)}{\inf_{\T\in\Delta(\gamma)}\mathcal{J}_{\mathrm{P}}(\T)} &=
1+o(1),\;\;\text{as}\;\;\gamma\to\infty,
\end{align*}
where $o(1)\to0$ as $\gamma\to\infty$.
\end{definition}

\begin{definition}[Second-Order Asymptotic Optimality]
A procedure $\T^*\in\Delta(\gamma)$ is said to be {\em second-order asymptotically} optimal in the class $\Delta(\gamma)$  if
\begin{align*}
\mathcal{J}_{\mathrm{P}}(\T^*)-\inf_{\T\in\Delta(\gamma)}\mathcal{J}_{\mathrm{P}}(\T)
&=O(1),\;\;\text{as}\;\;\gamma\to\infty,
\end{align*}
where $O(1)$ stays bounded as $\gamma\to\infty$.
\end{definition}

\begin{definition}[Third-Order Asymptotic Optimality]
A procedure $\T^*\in\Delta(\gamma)$ is said to be {\em third-order asymptotically} optimal in the class $\Delta(\gamma)$  if
\begin{align*}
\mathcal{J}_{\mathrm{P}}(\T^*)-
\inf_{\T\in\Delta(\gamma)}\mathcal{J}_{\mathrm{P}}(\T)
&=o(1),\;\;\text{as}\;\;\gamma\to\infty.
\end{align*}
\end{definition}

%-------------------------------------------------------------------------------------------------%
\subsection{The Shiryaev--Roberts--Pollak procedure}

The question of what procedure minimizes Pollak's measure of detection delay $\mathcal{J}_{\mathrm{P}}(T)$ is an open issue. As an attempt to resolve the issue,~\cite{Pollak:AS85} proposed to ``tweak'' the SR procedure~\eqref{eq:T-SR-def}. This led to the new procedure that we will refer to as the Shiryaev--Roberts--Pollak (SRP) procedure. To facilitate the presentation of the latter, we first explain the heuristics.

As known from the general decision theory (see, e.g.,~\citealp[Theorem~2.11.3]{Ferguson:Book67}), an $\mathcal{F}_n$-adapted stopping time $\T$ solves~\eqref{eq:Pollak-minimax-problem} if
\begin{inparaenum}[\itshape a\upshape)]
     \item $\T$ is an extended Bayes rule,
     \item it is an equalizer, and
     \item it satisfies the false alarm constraint with equality.
\end{inparaenum}
A procedure is said to be an equalizer if its conditional risk (which we measure through $\EV_{\nu}[\T-\nu|\T>\nu]$) is constant for all $\nu\ge 0$, that is, $\EV_{0}[\T]=\EV_{\nu}[\T-\nu|\T>\nu]$ for all $\nu\ge1$. Of the three conditions the one that requires $\T$ to be an equalizer poses the most challenge.~\Citet{Pollak:AS85} came up with an elegant solution.

It turns out that the sequence $\EV_{\nu}[\mathcal{S}_A-\nu|\mathcal{S}_A>\nu]$ indexed by $\nu$ eventually {\em stabilizes}, i.e., it remains the same for all sufficiently large $\nu$; see~\autoref{f:Fig1} below. This happens because the SR detection statistic enters the quasi-stationary mode, which means that the conditional distribution $\Pr_{\tinyinfty}(R_n\le x|\mathcal{S}_A>n)$ no longer changes with time. If one could get to the quasi-stationary mode immediately, then the resulting procedure would have the same expected conditional detection delay for  all $\nu \ge 0$, i.e., it would be the equalizer. Thus, Pollak's~\citeyearpar{Pollak:AS85} idea was to start the SR detection statistic $\{R_n\}_{n\ge0}$, defined in~\eqref{eq:Rn-SR-def}, not from zero ($R_0=0$), but from a random point $R_0=R_0^Q$, where $R_0^Q$ is sampled from the {\em quasi-stationary distribution} of the SR statistic under the hypothesis $\mathcal{H}_{\tinyinfty}$ (which is a Markov Harris-recurrent process under $\mathcal{H}_{\tinyinfty}$). Specifically, the quasi-stationary cdf $Q_A(x)$ is defined as
\begin{align}\label{eq:Q-SR-def}
Q_A(x)
&=
\lim_{n\to\infty}\Pr_{\tinyinfty}(R_n\le x|\mathcal{S}_A>n).
\end{align}

Therefore, the SRP procedure is defined as the stopping time
\begin{align}\label{eq:T-SRP-def}
\mathcal{S}_A^Q
&=
\inf\{n\ge1\colon R_n^Q\ge A\},
\end{align}
where $A>0$ is the detection threshold, and
\begin{align}\label{eq:R-SRP-def}
R_n^Q
&=
(1+R_{n-1}^Q)\LR_n, \quad n\ge 1, \quad R_0^Q\thicksim Q_A(x)
\end{align}
is the detection statistic.

We reiterate that, by design, the SRP procedure~\eqref{eq:T-SRP-def} and~\eqref{eq:R-SRP-def} is an equalizer: it delivers the same  conditional average detection delay for any change-point $\nu$, that is, $\EV_{0}[\mathcal{S}_A^Q]=\EV_{\nu}[\mathcal{S}_A^Q-\nu|\mathcal{S}_A^Q>\nu]$ for all $\nu\ge1$.

\Citet{Pollak:AS85} was able to demonstrate that the SRP procedure is third-order asymptotically optimal with respect to $\mathcal{J}_{\mathrm{P}}(\T)$. More specifically, the following is true.

\begin{theorem}[\citealp{Pollak:AS85}]
Let $\EV_{0}[(\log\LR_1)^+]<\infty$. Suppose that in the SRP procedure $\mathcal{S}_A^Q$ the detection threshold $A=A_{\gamma}$ is selected in such a way that $\EV_{\tinyinfty}[\mathcal{S}_{A_\gamma}^Q]=\gamma$. Then
\begin{align*}
\mathcal{J}_{\mathrm{P}}(\mathcal{S}_{A_\gamma}^Q)
&=
\inf_{\T\in\Delta(\gamma)}\mathcal{J}_{\mathrm{P}}(\T)+o(1),
\;\;\text{as}\;\;\gamma\to\infty.
\end{align*}
\end{theorem}
Recently,~\cite{Tartakovskyetal:TPA11} obtained the following asymptotic approximation for $\mathcal{J}_{\mathrm{P}}(\mathcal{S}_A^Q)$ under the second moment condition $\EV_{0}[\log\LR_1]^2<\infty$:
\begin{align*}
\EV_{0}[\mathcal{S}_{A}^Q]
&=
\frac{1}{I}(\log A+\varkappa - C_{\tinyinfty})+o(1),
\;\;\text{as}\;\;A\to\infty,
\end{align*}
where $\varkappa$ is the limiting average overshoot in the one-sided sequential test which is a subject of renewal theory (see, e.g.,~\citealp{Woodroofe:Book82}) and $C_{\tinyinfty}$ is a constant that can be computed numerically (e.g., by Monte Carlo simulations). Both $\varkappa$ and $C_{\tinyinfty}$ are formally defined in the next subsection, where we reiterate the exact result of~\cite{Tartakovskyetal:TPA11}.

Note that for sufficiently large $\gamma$,
\begin{align}\label{ARLSRP}
\EV_{\tinyinfty}[\mathcal{S}_A^Q]
&\approx
(A/\zeta)-\mu_Q,
\;\;\text{where}\;\;
\mu_Q=\int_0^A y\,dQ_A(y),
\end{align}
i.e., $\mu_Q$ is the mean of the quasi-stationary distribution, and $\zeta$ is a constant defined in~\eqref{kappazeta} below. This approximation can be obtained by first noticing that for a fixed $R_0^Q=r$ the process $R_n^Q-r-n$ is a zero-mean $\Pr_\infty$-martingale, and then applying optional sampling theorem to this martingale as well as a renewal theoretic argument  (cf.~\citealp{Tartakovskyetal:TPA11}).

%-------------------------------------------------------------------------------------------------%
\subsection{The Shiryaev--Roberts--$r$ procedure}\label{ss:SRr}

The third-order asymptotic optimality of the SRP procedure makes the latter practically appealing. On the flip size, the SRP rule requires the knowledge of the quasi-stationary distribution~\eqref{eq:Q-SR-def}. It is rare that this distribution can be expressed in a closed form; for examples where this is possible, see~\cite{Pollak:AS85},~\cite{Mevorach+Pollak:AJMMS91},~\cite{Polunchenko+Tartakovsky:AS10} and~\cite{Tartakovsky+Polunchenko:IWAP10}. As a result, the SRP procedure has not been used in practice.

To make the SRP procedure implementable,~\cite{Moustakidesetal-SS11} proposed a numerical framework. More importantly,~\cite{Moustakidesetal-SS11} offered numerical evidence that there exist procedures that are uniformly better than the SRP procedure. Specifically, they regard starting off the original SR procedure at a fixed (but specially designed) $R_0=r$, $0\le r<A$, and defining the stopping time with this new deterministic initialization. Because of the importance of the starting point, they dubbed their procedure the SR--$r$ procedure.

Formally, the SR--$r$ procedure is defined as the stopping time
\begin{align}\label{eq:T-SR-r-def}
\mathcal{S}_A^r
&=
\inf\{n\ge1\colon R_n^r\ge A\},
\end{align}
where $A>0$ is the detection threshold, and
\begin{align}\label{eq:R-SR-r-def}
R_n^r
&=
(1+R_{n-1}^r)\LR_n,
\;\; n\ge 1,\;\;\text{with}\;\; R_0^r=r\ge0
\end{align}
is the SR--$r$ detection statistic.

\Citet{Moustakidesetal-SS11} show numerically that for certain values of the starting point, $R_0^r=r$, apparently, $\EV_{\nu}[\mathcal{S}_{A_1}^r-\nu|\mathcal{S}_{A_1}^r>\nu]$ is strictly less than $\EV_{\nu}[\mathcal{S}_{A_2}^Q-\nu|\mathcal{S}_{A_2}^Q>\nu]$ for all $\nu\ge 0$, where $A_1$ and $A_2$ are such that $\EV_{\tinyinfty}[\mathcal{S}_{A_1}^r]=\EV_{\tinyinfty}[\mathcal{S}_{A_2}^Q]$ (although the maximal expected delay is only slightly smaller for $T_{A_1}^r$).
\begin{figure}
    \centering
    \includegraphics[width=0.95\textwidth]{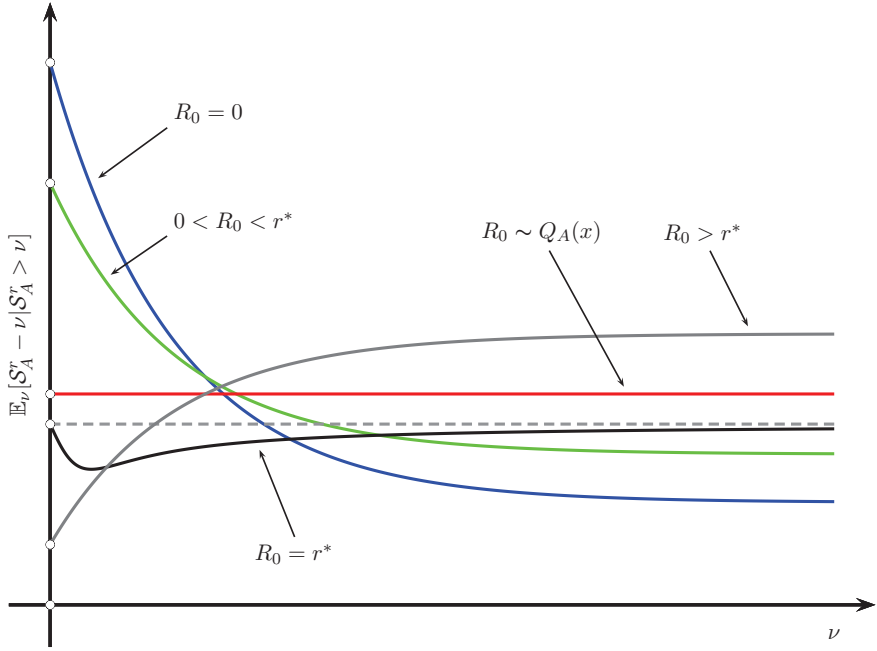}
    % figure caption to be below the figure
    \caption{Typical behavior of the conditional expected detection delay $\EV_{\nu}[\mathcal{S}_A^r-\nu|\mathcal{S}_A^r>\nu]$ of the SR--$r$ procedure as a function of the change-point $\nu$ for various initialization strategies.}
    \label{f:Fig1}
\end{figure}

It turns out that using the ideas of~\cite{Moustakidesetal-SS11} we are able to design the initialization point $r=r(\gamma)$ in the SR--$r$ procedure~\eqref{eq:T-SR-r-def} so that this procedure is also third-order asymptotically optimal. In this respect, the average delay to detection at infinity $\ADD_{\tinyinfty}(\mathcal{S}_A^r)=\lim_{\nu\to\infty}\EV_{\nu}[\mathcal{S}_A^r-\nu|\mathcal{S}_A^r>\nu]$ plays the critical role. To understand why, let us look at Figure~\ref{f:Fig1} which shows the average delay to detection  $\EV_{\nu}[\mathcal{S}_A^r-\nu|\mathcal{S}_A^r>\nu]$ vs. $\nu$ for several initialization values $R_0^r=r$. This figure was obtained using the integral equations and the numerical technique of~\cite{Moustakidesetal-SS11}. For $r=0$, this is the classical SR procedure (with $R_0=0$) whose average delay to detection  is monotonically decreasing to its minimum  (steady state value) that is attained at infinity. Note that in fact this steady state is attained for essentially finite values of the change-point $\nu$. It is seen that there exists a value $r=r_A$ that depends on the threshold $A$ for which the worst point $\nu$ is at infinity, i.e., $\mathcal{J}_{\mathrm{P}}(\mathcal{S}_A^{r_A})=\ADD_{\tinyinfty}(\mathcal{S}_A^{r_A})$. The value of $r_A$ is the minimal value for which this happens and it is also the value that delivers the minimum to the difference between $\mathcal{J}_{\mathrm{P}}(\mathcal{S}_A^{r}) $ and the lower bound for $\inf_{\T\in\Delta(\gamma)}\mathcal{J}_{\mathrm{P}}(\T)$ derived by~\cite{Moustakidesetal-SS11} and by~\cite{Polunchenko+Tartakovsky:AS10}. This is a very important observation, since it allows us to build a proof of asymptotic optimality based on an estimate of $\ADD_{\tinyinfty}(\mathcal{S}_A^{r})$. We also note that for the SR--$r$ procedure with initialization $r=r_A$ (pink line) the average detection delay at the beginning and at infinity are approximately equal, $\EV_0[\mathcal{S}_A^{r_A}]\approx\ADD_{\tinyinfty}(\mathcal{S}_A^{r_A})$. This allows us to conjecture that an ``optimal'' SR--$r$ is an equalizer at the beginning ($\nu=0$) and at sufficiently large values of $\nu$, so that initialization $r_A$ should be selected to achieve this property. The following theorem, whose proof can be found in \cite{Polunchenko+Tartakovsky:AS10}, shows that the lower bound for the ``minimax risk'' can be expressed via the integral average detection delay  of the SR--$r$ procedure, which partially explains the issue.

\begin{theorem}\label{Th1}
Let $\mathcal{S}_A^r$ be defined as in~\eqref{eq:T-SR-r-def} and~\eqref{eq:R-SR-r-def}, and let $A=A_\gamma$ be selected so that $\EV_{\tinyinfty}[\mathcal{S}_{A_\gamma}^r]=\gamma$. Then, for {\em every} $r\ge0$,
\begin{align}\label{IntADD}
\inf_{\T\in\Delta(\gamma)}\mathcal{J}_{\mathrm{P}}(\T)
&\ge
\frac{r\EV_{0}[\mathcal{S}_{A_\gamma}^r]+\sum_{\nu=0}^\infty \EV_{\nu}[(\mathcal{S}_{A_\gamma}^r-\nu)^+]}{r+\EV_{\tinyinfty}[\mathcal{S}_{A_\gamma}^r]}=
\mathcal{J}_{\mathrm{B}}(\mathcal{S}_{A_\gamma}^r).
\end{align}
\end{theorem}

Note that~\autoref{Th1} suggests that if $r$ can be chosen so that the SR--$r$ procedure is an equalizer (i.e., $\EV_{0}[\mathcal{S}_A^r]=\EV_\nu[\mathcal{S}_A^r-\nu|\mathcal{S}_A^r>\nu]$ for all $\nu \ge 0$), then it is {\em exactly} optimal. This is because the right-hand side in~\eqref{IntADD} is equal to $\EV_{0}[\mathcal{S}_A^r]$, which, in turn, is equal to $\sup_\nu \EV_\nu[\mathcal{S}_A^r-\nu|\mathcal{S}_A^r>\nu]=\mathcal{J}_{\mathrm{P}}(\mathcal{S}_A^r)$.
Therefore, we have the following corollary that will be used in Section~\ref{s:SR-r-opt} for proving that the SR--$r$ procedure with a specially designed $r=r_A$ is strictly optimal for two specific models.

\begin{corollary}\label{Cor:Cor1}
Let $A=A_\gamma$ be selected so that $\EV_{\tinyinfty}[\mathcal{S}_{A_\gamma}^r]=\gamma$. Assume that $r=r(\gamma)$ is chosen in such a way that the SR--$r$ procedure $\mathcal{S}_{A_\gamma}^{r(\gamma)}$ is an equalizer. Then it is strictly minimax in the class $\Delta(\gamma)$, i.e.,
\begin{align}\label{Minimax}
\inf_{\T\in\Delta(\gamma)}\mathcal{J}_{\mathrm{P}}(\T)
&=
\mathcal{J}_{\mathrm{P}}(\mathcal{S}_{A_\gamma}^{r(\gamma)}).
\end{align}
\end{corollary}

While the SR--$r$ is not strictly minimax in general, it is almost obvious that this procedure is almost minimax. In fact,~\cite{Moustakidesetal-SS11} conjecture that the SR--$r$ procedure is third-order asymptotically minimax and~\cite{Tartakovskyetal:TPA11} show that this conjecture is true. We will state the exact result after we introduce some additional notation.

Let $S_n=\log\LR_1+\cdots+\log\LR_n$ and, for $a \ge 0$, introduce the one-sided stopping time $\tau_a=\inf\{n\ge1\colon S_n \ge a\}$. Let $\kappa_a=S_{\tau_a}-a$ be an overshoot (excess over the level $a$ at stopping), and let
\begin{align}\label{kappazeta}
\varkappa&=\lim_{a\to\infty}\EV_{0}[\kappa_a], \quad \zeta= \lim_{a\to\infty}\EV_{0}\left[e^{-\kappa_a}\right].
\end{align}

The constants $\varkappa>0$ and $0<\zeta < 1$ depend on the model and can be computed numerically. Next, let $I=\EV_{0}[\log\LR_1]$ denote the Kullback--Leibler information number, and let $\tilde{V}_{\tinyinfty}=\sum_{j=1}^{\infty} e^{-S_j}$. Also, let $R_{\tinyinfty}$ be a random variable that has the $\Pr_{\tinyinfty}$-limiting (stationary) distribution of $R_n$, as $n \to\infty$, i.e., $Q_{\mathrm{ST}}(x)=\lim_{n\to\infty}\Pr_{\tinyinfty}(R_n\le x)=\Pr_{\tinyinfty}(R_{\tinyinfty}\le x)$. Let
\begin{align}\label{constant1}
C_{\tinyinfty}
&=
\EV[\log(1+R_{\tinyinfty}+\tilde{V}_{\tinyinfty})]
=
\int_0^\infty\int_0^\infty\log(1+x+y)\,dQ_{\mathrm{ST}}(x)\,d\tilde{Q}(y),
\end{align}
where $\tilde{Q}(y)=\Pr_{0}(\tilde{V}_{\tinyinfty}\le y)$.

\begin{theorem}[\citealp{Tartakovskyetal:TPA11}]
Let $\EV_{0}[\log\LR_1]^2<\infty$ and let $\log\LR_1$ be non-arithmetic. Then the following assertions hold.
\begin{enumerate}[\rm (i)]
    \item As $\gamma \to \infty$,
\begin{align}\label{LBasym}
\inf_{\T\in\Delta(\gamma)}\mathcal{J}_{\mathrm{P}}(\T)
&\ge
\frac{1}{I}[\log(\gamma\zeta)+\varkappa-C_{\tinyinfty}]+o(1).
\end{align}
    \item For any $r\ge0$,
\begin{align}\label{ADDinfty}
\ADD_{\tinyinfty}(\mathcal{S}_A^r)
&=
\EV_{0}[\mathcal{S}_A^{Q}]
=
\frac{1}{I}(\log A+\varkappa-C_{\tinyinfty})+o(1),
\;\;\text{as}\;\;A\to\infty.
\end{align}
\item Furthermore, if in the SR--$r$ procedure $A=A_\gamma=\gamma \zeta$ and the initialization point $r=o(\gamma)$ is selected so that $\mathcal{J}_{\mathrm{P}}(\mathcal{S}_A^r)=\ADD_{\tinyinfty}(\mathcal{S}_A^r)$, then $\EV_{\tinyinfty}[\mathcal{S}_A^r] = \gamma (1+o(1))$ and
\begin{align}\label{SADDSR-r}
\mathcal{J}_{\mathrm{P}}(\mathcal{S}_A^r)
&=
\frac{1}{I}[\log(\gamma \zeta)+\varkappa-C_{\tinyinfty}]+o(1)
\;\;\text{as}\;\;\gamma\to\infty.
\end{align}
Therefore, the SR--$r$ procedure is third-order asymptotically optimal:
\begin{align*}
\mathcal{J}_{\mathrm{P}}(\mathcal{S}_A^r)- \inf_{T\in\Delta(\gamma)}\mathcal{J}_{\mathrm{P}}(T)
&=o(1),
\;\;\text{as}\;\;\gamma\to\infty.
\end{align*}
\end{enumerate}
\end{theorem}

Also,
\begin{align}\label{SRrADD0}
\ADD_0(\mathcal{S}_A^r)
&=
\frac{1}{I}[\log A+\varkappa-C(r)]+o(1),
\;\;\text{as}\;\;
A\to\infty ,
\end{align}
where
\begin{align}\label{constCr}
C(r)
&=
\EV[\log(1+ r+\tilde{V}_{\tinyinfty})]=\int_0^\infty\log(1+r+y)\,d\tilde{Q}(y).
\end{align}
As we mentioned above, it is desirable to make the SR--$r$ procedure to look like equalizer by choosing the head start $r$, which can be achieved by equalizing
$\ADD_0$ and $\ADD_\infty$. Comparing~\eqref{ADDinfty} and~\eqref{SRrADD0} we see that this property approximately holds  when $r$ is selected from the equation $C(r^*)=C_\infty$. This shows that asymptotically as $\gamma\to\infty$ the ``optimal'' value $r^*$ is a fixed number that does not depend on $\gamma$. Clearly, this observation is important since it allows us to design the initialization point effectively and make the resulting procedure approximately optimal.

It is worth mentioning that for the conventional SR procedure that starts from zero
\begin{align*}
\mathcal{J}_{\mathrm{P}}(\mathcal{S}_A)
&=
\ADD_0(\mathcal{S}_A)=\frac{1}{I}[\log A+\varkappa-C(0)]+o(1),
\;\;\text{as}\;\;A\to\infty.
\end{align*}

Therefore, the SR procedure is only second-order asymptotically optimal. For sufficiently large $\gamma$, the difference between the supremum ADD-s of the SR procedure and the optimized SR--$r$ is given by $(C(0)-C_{\tinyinfty})/I$, which can be quite large if the Kullback--Leibler information number $I$ is small.

Note that similar to~\eqref{ARLSRP}, for sufficiently large $\gamma$,
\begin{align}\label{ARLSRr}
\EV_{\tinyinfty}[\mathcal{S}_A^r]&\approx(A/\zeta)- r.
\end{align}

\Citet{Polunchenko+Tartakovsky:AS10} and~\cite{Tartakovsky+Polunchenko:IWAP10} offer two scenarios where the SR--$r$ procedure is {\em strictly} minimax. Both are discussed (and extended) in~\autoref{s:SR-r-opt}. In addition, in~\autoref{sec:case-studies} we present an example where distributions $Q_{\mathrm{ST}}(x)$ and $\tilde{Q}(x)$ and the constants $\varkappa$, $\zeta$, $C_{\tinyinfty}$, and $C(r)$ can be computed analytically.

%-------------------------------------------------------------------------------------------------%
\section{Numerical performance evaluation}
\label{s:num-perf-eval} % NOT REFERENCED

Recall that each of the four approaches adduced above is characterized by its own optimality criterion. Together they bring about a variety of performance measures. Hence, to judge about the efficiency of a detection procedure (with respect to one performance measure or another), one has to be able to compute the procedure's corresponding operating characteristics (OC). In this section, we present integral equations for a multitude of OC-s that are of interest in various problem settings (Bayesian, minimax, etc.) and practical applications. Usually these equations cannot be solved analytically and numerical techniques are needed; cf.~\cite{Moustakidesetal-SS11,Tartakovskyetal-IWSM09}.

Consider a generic detection procedure described by the stopping time
\begin{align}\label{eq:generic-T-def}
\mathcal{T}_A^s
&=
\inf\{n\ge1\colon V_n^s\ge A\},
\end{align}
where $A>0$ is the detection threshold and $\{V_n^s\}_{n\ge0}$ is a Markov detection statistic that satisfies the recursive relation
\begin{align}\label{eq:generic-V-def}
V_n^s&=\xi(V_{n-1}^s)\LR_n,\;\; n\ge 1\;\;\text{with}\;\; V_{0}^s=s\ge0,
\end{align}
where $\xi(x)$ is a positive-valued function and $s$ is a fixed parameter referred to as the starting point or the ``head start''. Observe first that this detection scheme constitutes a rather broad class of detection schemes that includes, e.g., CUSUM, Shiryaev's procedure, SR--$r$, and EWMA (exponentially weighted moving average).  Indeed, for the Shiryaev procedure $\xi(V)=(1+V)/(1-p)$ and for the SR--$r$ procedure $\xi(V)=1+V$. (We will not consider other procedures such as CUSUM and EWMA in this paper).

Let $P_d^{\LR}(t)=\Pr_d(\LR_1\le t)$ denote the cdf of the LR $\LR_1$ under the measure $\Pr_d$, $d=\{0,\infty\}$. Define
\begin{align*}
\mathcal{K}_d(x,y)
&=
\frac{\partial}{\partial y}\Pr_d(V_{n+1}^s\le y|V_n^s=x)=\frac{\partial}{\partial y}P_d^{\LR}\left(\frac{y}{\xi(x)}\right),
\;\; d=\{0,\infty\},
\end{align*}
where $\{V_n^s\}_{n\ge0}$ is as in~\eqref{eq:generic-V-def}. Here and in the following we assume that $\LR_1$ is continuous under both hypotheses.

Let $\ell(x)=\EV_{\tinyinfty}[\mathcal{T}_A^x]$ and $\delta_{0}(x)=\EV_{0}[\mathcal{T}_A^x]$, where $\mathcal{T}_A^x$ is as in~\eqref{eq:generic-T-def}. Observe that $\ell(x)$ and  $\delta_{0}(x)$ are conditional expectations of the form $\EV_d[\,\cdot\,|V_{0}^x=x]$ (for $d=\infty$ and $d=0$ respectively), where $V_{0}^x=x$ is the starting point of the generic detection statistic~\eqref{eq:generic-V-def}. It is not difficult to see that $\ell(x)$ and $\delta_{0}(x)$ are governed by the equations
\begin{align}\label{eq:ARL-int-eqn}
\ell(x)
&=
1+\int_0^A\mathcal{K}_{\tinyinfty}(x,y)\,\ell(y)\,dy,
\end{align}
and
\begin{align}\label{eq:ADD0-int-eqn}
\delta_{0}(x)
&=
1+\int_0^A\mathcal{K}_{0}(x,y)\,\delta_{0}(y)\,dy,
\end{align}
respectively; cf.~\cite{Moustakidesetal-SS11}.

Next, for any $\nu\ge0$, let $\delta_{\nu}(x)=\EV_{\nu}[(\mathcal{T}_A^x-\nu)^+]$ and $\rho_{\nu}(x)=\Pr_{\tinyinfty}(\mathcal{T}_A^x>\nu)$. By~\cite{Moustakidesetal-SS11},
\begin{align}\label{eq:delta-rho-iter}
\delta_{\nu+1}(x)=\int_0^A\mathcal{K}_{\tinyinfty}(x,y)\,\delta_{\nu}(y)\,dy, \quad \rho_{\nu+1}(x)=\int_0^A\mathcal{K}_{\tinyinfty}(x,y)\,\rho_{\nu}(y)\,dy,
\end{align}
where $\delta_{0}(x)$ is governed by~\eqref{eq:ADD0-int-eqn} and $\rho_{0}(x) =1$ for all $x$, since $\Pr_{\tinyinfty}(\mathcal{T}_A^x>0)=1$. Consider now the conditional average delays to detection $\EV_{\nu}[\mathcal{T}_A^x-\nu|\mathcal{T}_A^x>\nu]=\EV_{\nu}[(\mathcal{T}_A^x-\nu)^+]/\Pr_{\nu}(\mathcal{T}_A^x>\nu)$, $\nu \ge 0$. Since $\Pr_{\nu}(\mathcal{T}_A^x>\nu)=\Pr_{\tinyinfty}(\mathcal{T}_A^x>\nu)$, we obtain
\begin{align*}
\EV_{\nu}[\mathcal{T}_A^x-\nu|\mathcal{T}_A^x>\nu]
&=
\frac{\EV_{\nu}[(\mathcal{T}_A^x-\nu)^+]}{\Pr_{\tinyinfty}(\mathcal{T}_A^x>\nu)} =\frac{\delta_\nu(x)}{\rho_\nu(x)}, \;\; \nu\ge0,
\end{align*}
where $\delta_{\nu}(x)$ and $\rho_{\nu}(x)$ are given by~\eqref{eq:delta-rho-iter}. Thus, the conditional average detection delays can be computed for any $\nu\ge0$, which allows one to evaluate $\sup_{\nu \ge 0} \EV_{\nu}[\mathcal{T}_A^x-\nu|\mathcal{T}_A^x>\nu]$.

Now, let $\psi(x)=\sum_{\nu=0}^{\infty}\EV_{\nu}[(\mathcal{T}_A^x-\nu)^+]=\sum_{\nu=0}^{\infty}\delta_{\nu}(x)$. By~\autoref{Th:SRoptimality},
\begin{align*}
\mathcal{J}_{\mathrm{ST}}(\mathcal{T}_A^x)
&=
\RIADD(\mathcal{T}_A^x)
=
\frac{\sum_{\nu=0}^{\infty}\EV_{\nu}[(\mathcal{T}_A^x-\nu)^+]}{\EV_{\tinyinfty}[\mathcal{T}_A^x)]}
=
\frac{\psi(x)}{\ell(x)},
\end{align*}
so that in order to compute the STADD $\mathcal{J}_{\mathrm{ST}}(\mathcal{T}_A^x)$ we have to be able to compute $\psi(x)$. As shown by~\cite{Moustakidesetal-SS11}, $\psi(x)$ is determined by the equation
\begin{align}\label{eq:IADD-int-eqn}
\psi(x)
&=
\delta_{0}(x)+\int_0^A\mathcal{K}_{\tinyinfty}(x,y)\,\psi(y)\,dy,
\end{align}
where $\delta_{0}(x)$ is governed by equation~\eqref{eq:ADD0-int-eqn}.

Note that the lower bound \eqref{IntADD} for the minimax risk given in Theorem~\ref{Th1} can be computed as
\begin{align*}
\mathcal{J}_{\mathrm{B}}(\mathcal{T}_A^x)
&=
\frac{x\delta_{0}(x)+\psi(x)}{x+\ell(x)},
\end{align*}
where $\ell(x)$, $\delta_{0}(x)$, and $\psi(x)$ are governed by equations~\eqref{eq:ARL-int-eqn},~\eqref{eq:ADD0-int-eqn}, and~\eqref{eq:IADD-int-eqn}.

The local conditional probabilities of false alarm $\Pr_{\tinyinfty}(\mathcal{T}_A^x\le k+m|\mathcal{T}_A^x> k)$, $k \ge 0$ inside a fixed ``window'' of size $m=1,2,\ldots$ can also be evaluated noting that $\Pr_{\tinyinfty}(\mathcal{T}_A^x\le k+m|\mathcal{T}_A^x> k)=1-\rho_{k+m}(x)/\rho_{k}(x)$,
where $\rho_{k}(x)$ are as in~\eqref{eq:delta-rho-iter}. Having $\Pr_{\tinyinfty}(\mathcal{T}_A^x\le k+m|\mathcal{T}_A^x>k)$ evaluated for sufficiently many $k$'s, one can easily find $\sup_{k}\Pr_{\tinyinfty}\mathcal{T}_A^x\le k+m|\mathcal{T}_A^x>k)$ for any fixed $m$.

The next step is to extend the obtained equations to the case when $\mathcal{T}_A^x$ is randomized similarly to the SRP procedure~\eqref{eq:T-SRP-def} and~\eqref{eq:R-SRP-def}. To this end, let $Q_A(y)=\lim_{n\to\infty}\Pr_{\tinyinfty}(V_n^s \le y|\mathcal{T}_A^s>n)$ be the quasi-stationary distribution. Note that this distribution does not depend on the starting point $V_0^s=s$ and  exists whenever the LR is continuous; cf.~\cite[Theorem~III.10.1]{Harris:Book63}.

It can be shown that the quasi-stationary pdf $q_A(x)=dQ_A(x)/dx$ satisfies the equation
\begin{align}\label{eq:QSD-int-eqn}
\lambda_A\,q_A(y)
&=
\int_0^Aq_A(x)\,\mathcal{K}_{\tinyinfty}(x,y)\,dx,\;\;\text{subject to}\;\;\int_0^A q_A(x)\,dx=1,
\end{align}
whence one can conclude that $q_A(x)$ is the {\em left dominant eigenvector} of the linear integral operator induced by the kernel $\mathcal{K}_{\tinyinfty}(x,y)$, and $\lambda_A\in(0,1)$ is the corresponding eigenvalue; cf.~\cite{Moustakidesetal-SS11} and~\cite{Pollak:AS85}. We also note that both $q_A(x)$ and $\lambda_A$ are {\em unique}.

Consider now $\mathcal{T}_A^Q$ defined as the above generic procedure $\mathcal{T}_A^x$ with the starting point being random and sampled from the quasi-stationary distribution. Specifically,
\begin{align}\label{eq:generic-T-rnd-def}
\mathcal{T}_A^Q
&=
\inf\{n\ge1\colon V_n^Q\ge A\},
\end{align}
where $A>0$ is the detection threshold, and $\{V_n^Q\}_{n\ge0}$ is a generic detection statistic computed recursively
\begin{align}\label{eq:generic-V-rnd-def}
V_n^Q
&=
\xi(V_{n-1}^Q)\LR_n,\;\; n\ge1\;\;\text{with}\;\; V_0^Q\thicksim Q_A.
\end{align}
We note that the SRP procedure is the special case of $\mathcal{T}_A^Q$ with $\xi(x)=1+x$.

Once $q_A(x)$ and $\lambda_A$ are available, one can compute the ARL to false alarm and the detection delay (which is independent from the change-point) for this randomized variant $\mathcal{T}_A^Q$ of the generic procedure $\mathcal{T}_A^x$. Indeed,
\begin{align*}
\EV_{\tinyinfty}[\mathcal{T}_A^Q]
&=
\int_0^A\ell(x)\,q_A(x)\,dx
=
\frac{1}{1-\lambda_A}
\;\;\text{and}\;\;
\EV_{0}[\mathcal{T}_A^Q]
=
\int_0^A\delta_{0}(x)\,q_A(x)\,dx .
\end{align*}
 To understand the second equality in the formula for $\EV_{\tinyinfty}[\mathcal{T}_A^Q]$, note that $\mathcal{T}_A^Q$ is $\Pr_{\tinyinfty}$-geometrically distributed with the ``probability of success'' $1-\lambda_A$. We also remark that, by design, the randomized variant $\mathcal{T}_A^Q$ of the generic procedure $\mathcal{T}_A^x$ is an equalizer, i.e., $\EV_{0}[\mathcal{T}_A^Q]=\EV_{\nu}[\mathcal{T}_A^Q-\nu|\mathcal{T}_A^Q>\nu]$ for all $\nu \ge 1$.

Finally, we present Bayesian operating characteristics --- the average detection delay
\begin{align*}
\ADD^\pi(\mathcal{T}_A^x)
&=
\frac{\sum_{k=0}^\infty\pi_k\EV_k(\mathcal{T}_A^x-k)^+}{1-\PFA^\pi(\mathcal{T}_A^x)}
\end{align*}
and the probability of false alarm, $\PFA^\pi(\mathcal{T}_A^x)= \sum_{k=1}^\infty\pi_k\Pr_{\tinyinfty}(\mathcal{T}_A^x \le k)$.  Assuming the geometric prior distribution~\eqref{eq:Shiryaev-geometric-prior-def}, we obtain
\begin{align*}
\begin{aligned}
\PFA^\pi(\mathcal{T}_A^x)
&=
(1-\pi)\left\{1-p\sum_{k=0}^\infty(1-p)^k\rho_k(x)\right\},
\\
\sum_{k=0}^\infty\pi_k\EV_k[(\mathcal{T}_A^x-k)^+]
&=
\pi\delta_0(r)+(1-\pi)p\sum_{k=0}^\infty(1-p)^k\delta_k(r)
\end{aligned}
\end{align*}
(cf.~\citealp{Tartakovsky+Moustakides:SA10}). Let $\psi_p(x)=\sum_{k=0}^\infty(1-p)^k\delta_k(x)$ and $\chi_p(x)=\sum_{k=0}^\infty(1-p)^k\rho_k(x)$.
Using the Markov property of the statistic $V_n^x$, it is readily seen that  $\psi_p(x)$ and $\chi_p(x)$ satisfy the following integral equations
\begin{align*}
\begin{aligned}
\psi_p(x)&=\delta_0(x)+(1-p)\int_0^A\mathcal{K}_{\tinyinfty}(x,y)\,\psi_p(y)\,dy,\\
\chi_p(x)&=1+(1-p)\int_0^A\mathcal{K}_{\tinyinfty}(x,y)\,\chi_p(y)\,dy .
\end{aligned}
\end{align*}

The PFA and ADD are then computed, respectively, as
\begin{align*}
\PFA^\pi(\mathcal{T}_A^x)&=(1-\pi)\left\{1-p\chi_p(x)\right\}
\;\;\text{and}\;\; \ADD^\pi(\mathcal{T}_A^x)=\frac{\pi\delta_0(x)+(1-\pi)p\psi_p(x)}{\pi+(1-\pi)p\chi_p(x)}.
\end{align*}

The above equations are Fredholm integral equations of the second kind. As a rule, such equations do not allow for an (exact) analytical solution. For a few exceptions from the rule see~\cite{Pollak:AS85},~\cite{Mevorach+Pollak:AJMMS91},~\cite{Polunchenko+Tartakovsky:AS10}, and~\cite{Tartakovsky+Polunchenko:IWAP10}. The results of the last two papers are summarized and extended in~\autoref{s:SR-r-opt}. Hence, a numerical technique may be in order. A simple numerical interpolation-projection type scheme has been suggested by~\cite{Moustakidesetal-SS11}. The scheme is effectively a piecewise collocation method with interpolating polynomials being of degree zero (constants). Using, e.g.,~\cite[Theorem~12.1.2]{Atkinson+Han:Book09} we can conclude that the corresponding rate of convergence is at worst linear.

The above performance evaluation methodology can now be applied to any particular scenario we may be interested in. A few such scenarios are worked out in Sections~\ref{s:SR-r-opt} and~\ref{sec:case-studies}.

%-------------------------------------------------------------------------------------------------%
\section{Exact optimality of the Shiryaev--Roberts--$r$ procedure}
\label{s:SR-r-opt} %

As we have pointed out earlier,  the question of  what solves Pollak's version of the minimax optimization problem~\eqref{eq:Pollak-minimax-problem}  has been open  since its inception in 1985. Because of the third-order asymptotic optimality and the fact that it is an equalizer it was conjectured that the SRP procedure might be the sought optimum. In this subsection, we suggest two counterexamples that disprove this conjecture. These examples show that a) the SRP procedure is not optimal, and b) that the SR--$r$ procedure is optimal. We stress that the SR--$r$ procedure is optimal in these examples, but not in general.

As a starting point, observe that equations~\eqref{eq:ARL-int-eqn},~\eqref{eq:ADD0-int-eqn},~\eqref{eq:IADD-int-eqn} and~\eqref{eq:QSD-int-eqn} are special cases of the more general equation
\begin{align}\label{eq:fredholm-int-eqn-2nd-kind}
u(x)
&=
v(x)+\int_0^A\mathcal{K}(x,y)\,u(y)\,dy,
\end{align}
where $v(x)$ is a given function, $u(x)$ is the sought (unknown) function, and $\mathcal{K}(x,y)$, which is called the {\em kernel} of this equation, is of the form
\begin{align*}
\mathcal{K}(x,y)
&=
\frac{\partial}{\partial y}P^{\LR}\left(\frac{y}{1+x}\right),
\end{align*}
with $P^{\LR}(x)$ being the cdf of the LR $\LR_n=g(X_n)/f(X_n)$.

To see that~\eqref{eq:fredholm-int-eqn-2nd-kind} is an ``umbrella'' equation for all equations we are interested in, note that to obtain equation~\eqref{eq:ARL-int-eqn},
which determines the ARL to false alarm, it suffices to take $v(x)=1$ for all $x$ and $\mathcal{K}(x,y)=\mathcal{K}_{\tinyinfty}(x,y)$. Likewise, equation~\eqref{eq:ADD0-int-eqn}, which governs the ADD at $\nu=0$, can be obtained from~\eqref{eq:fredholm-int-eqn-2nd-kind} by assuming $v(x)=1$ for all $x$ and $\mathcal{K}(x,y)=\mathcal{K}_{0}(x,y)$. By a similar argument, one can also verify that equations~\eqref{eq:IADD-int-eqn} and~\eqref{eq:QSD-int-eqn} are instances of~\eqref{eq:fredholm-int-eqn-2nd-kind}. Thus, if one is able to solve equation~\eqref{eq:fredholm-int-eqn-2nd-kind}, one is also able to solve any of the equations of interest.

Suppose now that we have a change-point scenario for which the cdf $P^{\LR}(t)$ is such that
\begin{align*}
P^{\LR}\left(\frac{y}{1+x}\right)
&=
\mathcal{X}(x)\,\mathcal{Y}(y)
\end{align*}
for some sufficiently smooth functions $\mathcal{X}(x)$ and $\mathcal{Y}(y)$. In this case, the kernel $\mathcal{K}(x,y)$ is separable, i.e.,
\begin{align*}
\mathcal{K}(x,y)
&=
\frac{\partial}{\partial y}P^{\LR}\left(\frac{y}{1+x}\right)
=
\mathcal{X}(x)\frac{d}{dy}\mathcal{Y}(y)
=
\mathcal{X}(x)\,\mathcal{Y}'(y),
\end{align*}
so that the variables $x$ and $y$ are separated.

If the kernel is separable and the interval of integration has constant limits, the above equation can be solved analytically, and the solution is $u(x)
=v(x)+M\mathcal{X}(x)$, where
\begin{align*}
M
&=
\left(\,\int_0^A v(t)\,\mathcal{Y}'(t)\,dt\right)\left/\left(1-\int_0^A\mathcal{X}(t)\,\mathcal{Y}'(t)\,dt\right)\right.,
\end{align*}
which is a function of $A$ only.

More important is the fact that in this case
\begin{align*}
\Pr(R_1^r\le y|\mathcal{S}_A^r>1)
&=
\Pr(R_1^r\le y|R_1^r<A,R_0^r=r)
=
\mathcal{Y}(y)/\mathcal{Y}(A),
\end{align*}
i.e., $\Pr(R_1^r\le y|\mathcal{S}_A^r>1)$  does not depend on the starting point $R_0^r=r$. This means that the quasi-stationary distribution ``kicks in'' as early as the first observation becomes available. As a result, the SR--$r$ procedure is an equalizer for $\nu \ge 1$, and the only ``degree of freedom'' is $\nu=0$. If one now designs the starting point $R_0^r$ so as to equate the performance of the SR--$r$ procedure at $\nu=0$ to that at $\nu \ge 1$, then the SR--$r$ procedure will be an equalizer for all $\nu \ge 0$.  Therefore, by~\autoref{Cor:Cor1}, in this case it is minimax. Note that this equalizer is different from the SRP rule, which is also an equalizer. The SR--$r$ is an equalizer and minimax not in general but only in this particular case, i.e., in the case when the kernel is separable. Thus, we now have to find examples where this is true. This will prove that the SRP procedure is not strictly minimax in general.

Suppose now that the observations' distribution is $\mathsf{uniform}(0,1)$ pre-change and $\mathsf{beta}(2,1)$ post-change, that is, $f(x)=\indicator{0<x<1}$ and $ g(x)=2x\indicator{0<x<1}$. The LR is $\LR_n=2X_n\indicator{0<X_n<1}$; observe that $\LR_n\in(0,2)$, since $X_n\in(0,1)$. Hence,
\begin{align*}
P_{\tinyinfty}^{\LR}(t)
&=
\begin{cases}
1,&\text{if}\,\, t \ge 2;\\
t/2,&\text{if}\,\,0\le t<2;\\
0,&\text{otherwise,}
\end{cases}
\quad\text{and}\quad
P_{0}^{\LR}(t)
=
\begin{cases}
1,&\text{if}\,\, t \ge 2;\\
(t/2)^2,&\text{if}\,\,0\le t<2;\\
0,&\text{otherwise.}
\end{cases}
\end{align*}
It is apparent that both these distributions are monomial and therefore separable. As a result, one can compute the required operating characteristics of any SR-type procedure analytically. This was done by~\cite{Tartakovsky+Polunchenko:IWAP10}. Another example, where $f(x)=e^{-x}\indicator{x\ge0}$ and $g(x)=2e^{-2x}\indicator{x\ge0}$, was considered by~\cite{Polunchenko+Tartakovsky:AS10}. Although this model may seem very different from the $\mathsf{uniform}(0,1)$-to-$\mathsf{beta}(2,1)$ model, it has exactly the same distributions $P_d^{\LR}(t)$, $d=\{0,\infty\}$. Both papers established the following theorem the proof of which can be found in~\cite{Polunchenko+Tartakovsky:AS10}.

\begin{theorem}[\citealp{Polunchenko+Tartakovsky:AS10}]
Let $\bar{\gamma}=1/(1-0.5\log3)\approx2.2$.
\begin{enumerate}[\rm (i)]
    \item If the starting point $r$ in the SR--$r$ procedure is chosen as $r_A=\sqrt{1+A}-1$ and the detection threshold $A=A_{\gamma}$ is set to the solution of the transcendental equation
\begin{align*}
A+(\gamma-1)\sqrt{1+A}\log(1+A)-2(\gamma-1)\sqrt{1+A}
&=
0,
\end{align*}
then, for every $\gamma\in(1,\bar{\gamma})$, the ARL to false alarm $\EV_{\tinyinfty}[\mathcal{S}_A^r]$ is exactly $\gamma$ and the SR--$r$ procedure is strictly minimax. That is,
\begin{align*}
\mathcal{J}_{\mathrm{P}}(\mathcal{S}_A^r)
&=\inf_{\T\in\Delta(\gamma)}\mathcal{J}_{\mathrm{P}}(\T)
\;\;\text{for every}\;\;\gamma\in(1,\bar{\gamma}).
\end{align*}
     \item If the detection threshold in the SRP procedure $\mathcal{S}_B^Q$ is set to $B=B_{\gamma}=\exp\{2(1-1/\gamma)\}-1$, then the ARL to false alarm $\EV_{\tinyinfty}[\mathcal{S}_B^Q]$ is exactly $\gamma$ and $\mathcal{J}_{\mathrm{P}}(\mathcal{S}_B^Q)$ is strictly greater than $\mathcal{J}_{\mathrm{P}}(\mathcal{S}_A^r)$ for every $\gamma\in(1,\bar{\gamma})$. Hence, the SRP procedure is suboptimal.
\end{enumerate}
\end{theorem}

\begin{figure}[h]
    \centering
    \includegraphics[width=0.85\textwidth]{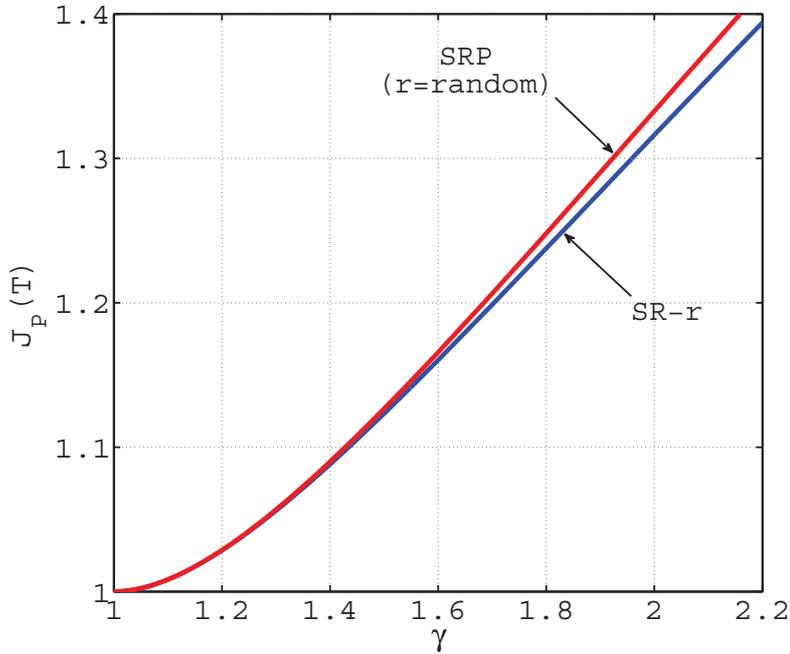}
    % figure caption to be below the figure
    \caption{Performance of the SRP procedure vs. that of the SR--$r$ rule for the $\textsf{uniform}(0,1)$-to-$\textsf{beta}(2,1)$ model. The detection threshold for either procedure is between $0$ and $2$.}
    \label{fig:SRP-SR-r_counterexample}
\end{figure}

This theorem is illustrated in~\autoref{fig:SRP-SR-r_counterexample}. Note that the curves in the picture are {\em exact}. We stress again that while the SR--$r$ procedure is exactly minimax in this example, it is still an open question what minimizes Pollak's $\mathcal{J}_{\mathrm{P}}(\T)$ in general. We conjecture that, in the general case, the optimal procedure is based on the deterministic initialization $\{r_n\}_{n\ge 1}$ that depends on time, in which case a threshold $A=A_n$ may also be a function of time. We also note that even if this fact is proved rigorously, finding the sequences $\{r_n(\gamma)\}$ and $\{A_n(\gamma)\}$ (in every particular case) is an extremely difficult problem. Solving this problem may not be worth trying since the difference between the lower bound~\eqref{IntADD} and the supremum ADD is usually small, at least for a moderate (and of course low) false alarm rate. See~\cite{Moustakidesetal-SS11} and \autoref{fig:exp-exp-scenario:STADD-vs-ARL} below.

We conclude this subsection with a remark concerning exact optimality of the SR--$r$ procedure in the class $\Delta_{\alpha}^m=\{\T\colon \sup_{k}\Pr_{\tinyinfty}(k<\T\le k+m|\T>k)\le\alpha\}$, where $\alpha\in(0,1)$ and $m\ge1$. We first discussed this class in~\autoref{sec:opt-criteria-overview:minimax}, where we mentioned that it is ``stronger'' than the class $\Delta(\gamma)$, i.e., in general $\Delta(\gamma)$ contains $\Delta_{\alpha}^m$. It can be easily verified that the $\Pr_{\tinyinfty}$-distribution of the SR--$r$ stopping time $\mathcal{S}_A^r$ is zero-modified geometric:
\begin{empheq}[%
    left={%
        \Pr_{\tinyinfty}(k<\mathcal{S}_A^r\le k+m|\mathcal{S}_A^r>k)=1-%
    \empheqlbrace}]{align*}
&\left[\frac{1}{2}\log(1+A)\right]^m\quad\text{for $k \ge 1$;}\\
&\frac{A}{2(1+r)}\left[\frac{1}{2}\log(1+A)\right]^{m-1}\quad\text{for $k=0$},
\end{empheq}
where $m\ge1$. Thus, there is a one-to-one correspondence between the classes $\Delta_{\alpha}^m$ and $\Delta(\gamma)$. As a result, the SR--$r$ procedure is minimax in the class $\Delta_{\alpha}^m$ as well. The same is true for the exponential model considered by~\cite{Polunchenko+Tartakovsky:AS10}. We believe that this is the first exact optimality result in the class $\Delta_{\alpha}^m$.

%-------------------------------------------------------------------------------------------------%
\section{Case studies}
\label{sec:case-studies} %

This section dissects two specific cases of the iid model~\eqref{iidmodel}  to illustrate the performance margin between the SR$-r$ and SRP procedures $\mathcal{S}_A^r$ and $\mathcal{S}_A^Q$, defined in~\eqref{eq:T-SR-r-def}, \eqref{eq:R-SR-r-def} and \eqref{eq:T-SRP-def}, \eqref{eq:R-SRP-def}, respectively.

%-------------------------------------------------------------------------------------------------%
\subsection{Example 1: A beta-to-beta model}

Suppose the pre- and post-change densities are, respectively,
\begin{align*}
f(x)
&=
\frac{x^{\delta-1}(1-x)^{\delta}}{\mathtt{B}(\delta,\delta+1)} \indicator{0<x<1}
\;\;\text{and}\;\;
g(x)
=
\frac{x^{\delta}(1-x)^{\delta-1}}{\mathtt{B}(\delta+1,\delta)}
\indicator{0<x<1},
\end{align*}
where $\delta>0$ is a given constant and $\mathtt{B}(\cdot,\cdot)$ is the Beta function. That is, the observations $X_n$, $n\ge1$ are iid $\mathtt{beta}(\delta,\delta+1)$-distributed pre-change and iid $\mathtt{beta}(\delta+1,\delta)$-distributed post-change. This  $\mathtt{beta}(\delta,\delta+1)$-to-$\mathtt{beta}(\delta+1,\delta)$ model is of interest in the context of studying the accuracy of the asymptotic expansions for the performance of the two competing SR-type procedures. Specifically, recall that for sufficiently large detection thresholds,
\begin{align*}
\EV_{\tinyinfty}[\mathcal{S}_A^Q]
&\approx A/\zeta-\mu_Q
\;\;\text{and}\;\;
\ADD_{\nu}(\mathcal{S}_A^Q)
\approx
\frac{1}{I}(\log A+\varkappa-C_{\tinyinfty})
\;\;\text{for all $\nu\ge0$,}
\\
\EV_{\tinyinfty}[\mathcal{S}_A^r]
&\approx
A/\zeta-r
\;\;\text{and}\;\;
\ADD_{\tinyinfty}(\mathcal{S}_A^r)
\approx
\frac{1}{I}(\log A+\varkappa-C_{\tinyinfty}),
\end{align*}
where $I=\EV_0[\log\LR_1]$ is the Kullback--Leibler information number, $\zeta$ and $\varkappa$ are defined in~\eqref{kappazeta}, $\mu_Q$ is the mean of the quasi-stationary distribution, and the constant $C_\infty$ is defined in~\eqref{constant1}.

For the $\mathtt{beta}(\delta,\delta+1)$-to-$\mathtt{beta}(\delta+1,\delta)$ model, $C_{\tinyinfty}$, $\varkappa$, $\zeta$ and $I$ are all computable {\em analytically} for any $\delta>0$. This is of much aid in the context of testing the accuracy of the asymptotic approximations. Specifically, we first present the exact, explicit formulas for each of the needed quantities, assuming arbitrary $\delta>0$. We then evaluate the performance of the procedures of interest using the methodology of~\autoref{s:num-perf-eval} and compare the obtained performance against that predicted by the asymptotic approximations.

Observe that $\LR_n=X_n/(1-X_n)$ for any $\delta>0$, whence one can readily deduce $P_{d}^{\LR}(t)=\Pr_{d}(\LR_1\le t)$, $d=\{0,\infty\}$. Specifically, the densities $p_d^{\LR}(t)=dP_{d}^{\LR}(t)/dt$, $d=\{0,\infty\}$, can be seen to be
\begin{align}\label{eq:beta2beta-LR-pdf-s}
p_{\tinyinfty}^{\LR}(t)
&=
\frac{t^{\delta-1}(1+t)^{-2\delta-1}}{\mathtt{B}(\delta,\delta+1)} \indicator{t>0}
\;\;\text{and}\;\;
p_{0}^{\LR}(t)
=
\frac{t^{\delta}(1+t)^{-2\delta-1}}{\mathtt{B}(\delta+1,\delta)} \indicator{t>0},
\end{align}
i.e., under either measure $\Pr_d$, $d=\{0,\infty\}$, the LR's distribution is Beta of type II (also known as the Beta prime distribution); the parameters are $\delta$ and $\delta+1$ under measure $\Pr_{\tinyinfty}$, and $\delta+1$ and $\delta$ under measure $\Pr_0$. The fact that $p_{\tinyinfty}^{\LR}(t)$ and $p_{0}^{\LR}(t)$ are both Beta prime with ``mirrored'' parameters suggests a certain symmetry embedded in the $\mathtt{beta}(\delta,\delta+1)$-to-$\mathtt{beta}(\delta+1,\delta)$ model. Specifically, consider the ``dual'' $\mathtt{beta}(\delta+1,\delta)$-to-$\mathtt{beta}(\delta,\delta+1)$ model. That is, suppose the pre- and post-change distributions -- $f(x)$ and $g(x)$ -- are swapped so that the former is not $\mathtt{beta}(\delta,\delta+1)$, but $\mathtt{beta}(\delta+1,\delta)$, and the latter is not $\mathtt{beta}(\delta+1,\delta)$, but $\mathtt{beta}(\delta,\delta+1)$. A special case of this swapped model (with $\delta=1$) was considered by~\cite{Tartakovskyetal:TPA11}. It can be shown, exploiting properties of the Beta and Beta prime distributions, that for the swapped $\mathtt{beta}(\delta+1,\delta)$-to-$\mathtt{beta}(\delta,\delta+1)$ model, the densities $p_d^{\LR}(t)$, $d=\{0,\infty\}$, are exactly the same as those we just derived for the original $\mathtt{beta}(\delta,\delta+1)$-to-$\mathtt{beta}(\delta+1,\delta)$ model; see~\eqref{eq:beta2beta-LR-pdf-s}. Put otherwise, the $\mathtt{beta}(\delta,\delta+1)$-to-$\mathtt{beta}(\delta+1,\delta)$ model and the $\mathtt{beta}(\delta+1,\delta)$-to-$\mathtt{beta}(\delta,\delta+1)$ model are statistically indistinguishable, for any $\delta>0$. This symmetry entails a few consequences to be demonstrated next.

Consider the stationary distribution $Q_{\mathrm{ST}}(x)=\lim_{n\to\infty}\Pr_{\tinyinfty}(R_n\le x)$ of the SR statistic $\{R_n\}_{n\ge0}$. The quasi-stationary pdf  $q_{\mathrm{ST}}(x)=dQ_{\mathrm{ST}}(x)/dx$ is governed by the equation
\begin{align*}
q_{\mathrm{ST}}(x)
&=
\int_0^\infty \frac{\partial}{\partial x}P_{\tinyinfty}^{\LR}\left(\frac{x}{1+y}\right) q_{\mathrm{ST}}(y)\, dy,
\end{align*}
which can be derived from equation \eqref{eq:QSD-int-eqn} for the quasi-stationary pdf, $q_A(x)$, by letting $A\to\infty$ and noticing that $\lim_{A\to\infty}\lambda_A=1$ and $\lim_{A\to\infty}q_A(x)=q_{\mathrm{ST}}(x)$ (cf.~\citealp{Pollak+Siegmund:JAP86}). Using \eqref{eq:beta2beta-LR-pdf-s}, we obtain \begin{align*}
q_{\mathrm{ST}}(x)
&=
\frac{x^{\delta-1}}{\mathtt{B}(\delta+1,\delta)}\int_0^\infty\frac{(1+y)^{\delta+1}}{(1+x+y)^{1+2\delta}} q_{\mathrm{ST}}(y) \, d y,
\end{align*}
and the (exact) solution is
\begin{align*}
q_{\mathrm{ST}}(x)
&=
\frac{x^{\delta-1}(1+x)^{-1-\delta}}{\mathtt{B}(\delta,1)}\indicator{x>0}
=
\delta x^{\delta-1}(1+x)^{-1-\delta}\indicator{x>0},
\end{align*}
which is the pdf of a Beta prime distribution with parameters $\delta$ and $1$. Note that $q_{\mathrm{ST}}(x)\thicksim x^{-2}$ as $x\to\infty$, which agrees with~\cite{Kesten:AM73}.

Next, it can be shown that the pdf  $\tilde{q}(x)=d\tilde{Q}(x)/dx$ of distribution $\tilde{Q}(x)=\Pr_{0}(\tilde{V}_{\tinyinfty}\le x)$ is governed by the equation
\begin{align*}
\tilde{q}(x)
&=
-\int_0^\infty \frac{\partial}{\partial x}P_{0}^{\LR}\left(\frac{1+y}{x}\right) \tilde{q}(y) \, dy,
\end{align*}
which can be established in a manner similar to that used to derive the above equation for $q_{\mathrm{ST}}(x)$. However, due to the symmetry of the model one can immediately conclude that $\tilde{q}(x)\equiv q_{\mathrm{ST}}(x)$, so that
\begin{align}\label{eq:q-tilde-pdf}
\tilde{q}(x)
&=
q_{\mathrm{ST}}(x)
=
\delta x^{\delta-1}(1+x)^{-1-\delta}\indicator{x>0} .
\end{align}
 We now can find
\begin{align*}
C_{\tinyinfty}
&=
\delta\Psi_1(\delta)+\Psi_0(\delta)-\Psi_0(1),
\end{align*}
where $\Psi_n(x)=d^{n+1}\log\Gamma(x)/dx^{n+1}$ ($n\ge0$) is the polygamma function and $\Gamma(x)$ is the Gamma function; also note that $\Psi_0(1)=-0.577\ldots$ is the negative Euler's constant.

To find $\zeta$ and $\varkappa$, we use the formulas
\begin{align*}
\zeta
&=
\frac{1}{I}\exp\left\{-\sum_{k=1}^\infty\frac{1}{k}\bigl[\Pr_{\tinyinfty}(S_k>0)+\Pr_0(S_k\le 0)\bigr]\right\},
\\
\varkappa &=
\frac{\EV_0[Z_1^2]}{2I}-\sum_{k=1}^\infty\frac{1}{k}\EV_0[S_k^-],
\end{align*}
where $x^-=-\min(0,x)$; cf., e.g.,~\cite[Chapters~2~\&~3]{Woodroofe:Book82} and~\cite[Chapter~VIII]{Siegmund:Book85}. Using the work of~\cite{Springer+Thompson:JAM70}, after certain manipulations we obtain that
\begin{align*}
\Pr_{0}(S_k\le0)
&=
\frac{1}{\Gamma^k(\delta)\Gamma^k(\delta+1)}G_{k+1,k+1}^{k+1,k}\left(1\left|
\begin{array}{c}
\overbrace{-\delta,\ldots,-\delta}^\textrm{$k$ times},1 \\ 0,\underbrace{\delta,\ldots,\delta}_\textrm{$k$ times}
\end{array}\right.\right),\;\; k\ge1,
\end{align*}
where $G_{\cdot,\cdot}^{\cdot,\cdot}(\cdot|\cdot)$ is the Meijer G-function. Note that due to the symmetry of the $\mathtt{beta}(\delta,\delta+1)$-to-$\mathtt{beta}(\delta+1,\delta)$ model,  $\Pr_{\tinyinfty}(S_k>0)=\Pr_0(S_k\le 0)$ for all $k\ge1$. Hence,
\begin{align*}
\zeta
&=
\delta\exp\left\{-2\sum_{k=1}^\infty\frac{1}{k\Gamma^k(\delta)\Gamma^k(\delta+1)}G_{k+1,k+1}^{k+1,k}\left(1\left|
\begin{array}{c}
-\delta,\ldots,-\delta,1 \\ 0,\delta,\ldots,\delta
\end{array}\right.\right)\right\},
\end{align*}
which can be evaluated numerically for any $\delta>0$ and with any desired accuracy.

Next, it can be shown that $\EV_0[Z_1^2]=2\Psi_1(\delta)$ and
\begin{align*}
\EV_0[S_k^-]
&=
\frac{1}{\Gamma^k(\delta)\Gamma^k(\delta+1)}G_{k+2,k+2}^{k+2,k}\left(1\left|
\begin{array}{c}
-\delta,\ldots,-\delta, 1, 1\\ 0, 0,\delta,\ldots,\delta
\end{array}\right.\right),\;\; k\ge1,
\end{align*}
whence
\begin{align*}
\varkappa
&=
\delta\Psi_1(\delta)
-\sum_{k=0}^\infty\frac{1}{k\Gamma^k(\delta)\Gamma^k(\delta+1)}G_{k+2,k+2}^{k+2,k}\left(1\left|
\begin{array}{c}
-\delta,\ldots,-\delta, 1, 1\\ 0, 0,\delta,\ldots,\delta
\end{array}\right.\right).
\end{align*}

Consider now starting the SR--$r$ procedure off the point $R_0^r=r^*$ for which $\ADD_0(\mathcal{S}_A^{r^*})$ and $\ADD_{\tinyinfty}(\mathcal{S}_A^{r^*})$ are the same (at least approximately). This idea was first brought up in~\autoref{ss:SRr}; recall~\autoref{f:Fig1}. By~\eqref{ADDinfty} and~\eqref{SRrADD0}, when the ARL to false alarm is sufficiently large,
\begin{align}\label{approxADD}
\ADD_{\tinyinfty}(\mathcal{S}_A^r)
&\approx
\frac{1}{I}(\log A+\varkappa-C_{\tinyinfty})
\;\;\text{and}\;\;
\ADD_{0}(\mathcal{S}_A^r)
\approx
\frac{1}{I}[\log A+\varkappa-C(r)],
\end{align}
where $C(r)=\EV[\log(1+r+\tilde{V}_{\tinyinfty})]$ (see \eqref{constCr}). Hence, equating $\ADD_0(\mathcal{S}_A^r)$ and $\ADD_{\tinyinfty}(\mathcal{S}_A^r)$ is equivalent to requiring $C_{\tinyinfty}=C(r)$, and setting $R_0^r$ to $r$ that solves the equation $C_{\tinyinfty}=C(r)$ results in the desired effect of $\ADD_0(\mathcal{S}_A^r)\approx\ADD_{\tinyinfty}(\mathcal{S}_A^r)$ (asymptotically). Since $C_{\tinyinfty}$ is already computed, it is left to find $C(r)$. To this end, using~\eqref{eq:q-tilde-pdf}, we obtain
\begin{align*}
C(r)
&=
\Phi\left(\frac{r}{1+r},1,\delta\right)+\Psi_0(\delta)-\Psi_0(1),
\end{align*}
where $\Phi(\cdot,\cdot,\cdot)$ is the Lerch transcendent. Hence, the equation $C_{\tinyinfty}=C(r)$, where $r$ is the unknown, reduces to
\begin{align*}
\Phi\left(\frac{r}{1+r},1,\delta\right)
&=
\delta\Psi_1(\delta),\;\; r\ge0,
\end{align*}
which can be solved numerically for any desired $\delta>0$ and with any pleased precision.

We are now in a position to perform particular computations.
To remind, we would like to test the accuracy of the asymptotic approximations~\eqref{ADDinfty} and~\eqref{SADDSR-r}. Clearly, the accuracy is the better, the higher the desired level of the ARL to false alarm $\EV_{\tinyinfty}[\T]=\gamma$. First, we intend to try a relatively small value of $\gamma=10^2$, which corresponds to practically high chances of sounding a false alarm. We do not expect the asymptotics to kick in for $\gamma$ lower than a few hundreds. Suppose that $\delta=1$. For this choice of $\delta$ we have: $C_{\tinyinfty}=\pi^2/6\approx1.64$, $I=1$, $\zeta\approx0.425$, $\varkappa\approx1.25$, and $r^*\approx2$ (so that $C(r^*)=C_{\tinyinfty}\approx1.64$).

The first step is to set thresholds to guarantee the given ARL to false alarm $\gamma$. For the SR--$r$ procedure,  the detection threshold, $A$, should be set to the solution of the equation $\gamma=A/\zeta-r$,  which follows  from the corresponding asymptotics for $\EV_{\tinyinfty}[\mathcal{S}_A^r]$. Since in our case $r=r^*\approx2$ and $\zeta\approx0.425$, we find that $A$ must be set to about $43$. The actual (evaluated numerically with very high accuracy) ARL to false alarm with this $A$ is $100.1$. Hence, the approximation $\EV_{\tinyinfty}[\mathcal{S}_A^r]\approx A/\zeta-r$ is very accurate, even when $\gamma=10^2$, which is equivalent to a relatively high risk of raising a false alarm. For the SRP procedure to have $\EV_{\tinyinfty}[\mathcal{S}_A^Q]=10^2$ the detection threshold, $A$, should be set to $43$ as well; the actual ARL to false alarm for this choice of $A$ is $99.6$, and the mean, $\mu_Q$, of the quasi-stationary distribution is around $2.6$. Hence, the approximation $\EV_{\tinyinfty}[\mathcal{S}_A^Q]\approx A/\zeta-\mu_Q$ is also very accurate.

\begin{figure}
    \centering
    \includegraphics[width=0.95\textwidth]{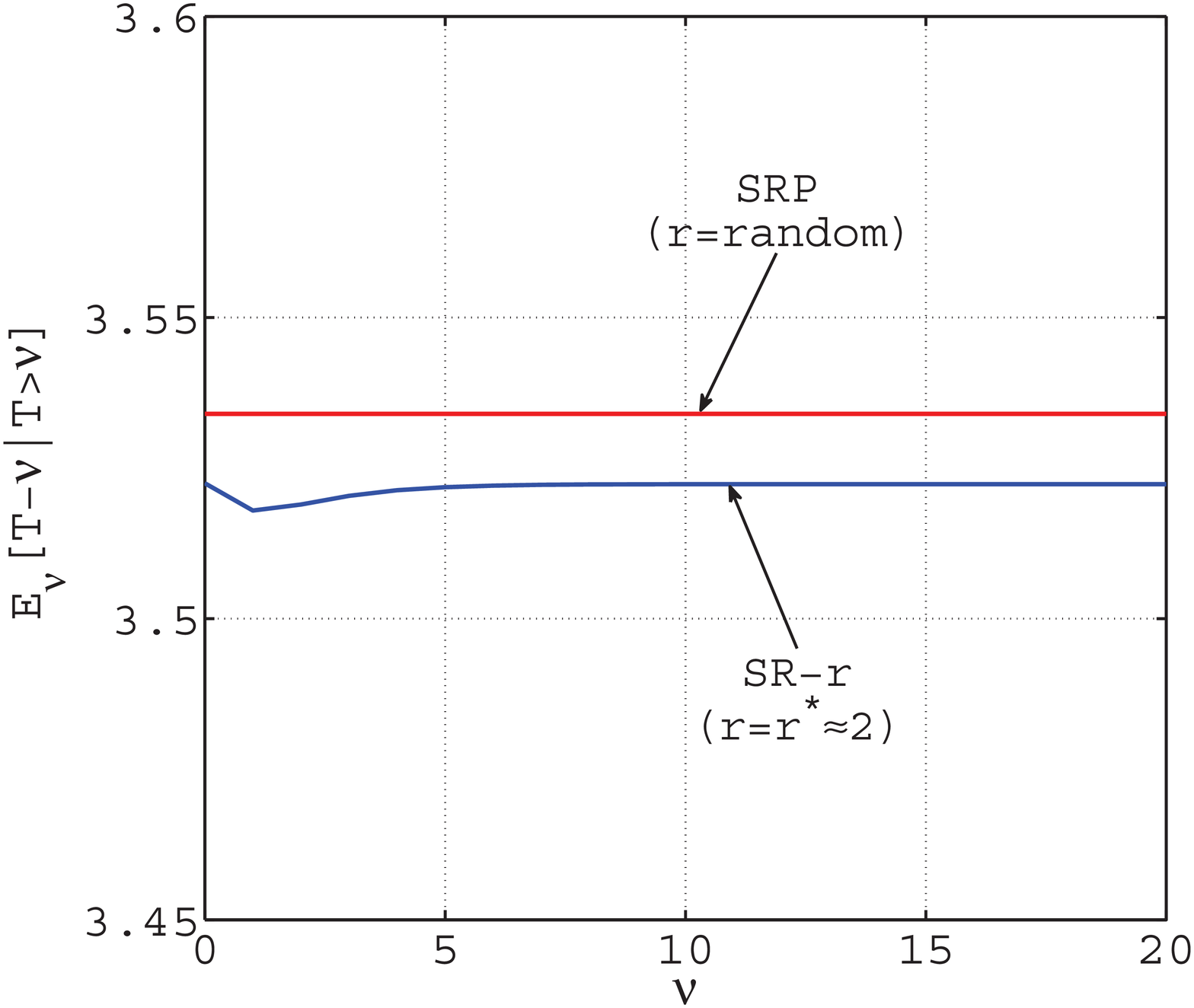}
    % figure caption to be below the figure
    \caption{Conditional average detection delay $\EV_{\nu}[\mathcal{S}_A^r-\nu|\mathcal{S}_A^r>\nu]$ vs. change-point $\nu$ for the SRP procedure and for the SR--$r$ procedure with $r=r^{*}\approx2$ for the $\mathsf{beta}(\delta,\delta+1)$-to-$\mathsf{beta}(\delta+1,\delta)$ model with $\delta=1$. The ARL to false alarm $\EV_{\tinyinfty}[\T]=\gamma$ is approximately $100$ for each procedure.}
    \label{fig:SRP_SR_r_ADDk_vs_k_ARL_100}
\end{figure}

We now proceed to examining $\ADD_{\nu}(\T)=\EV_\nu[\T-\nu|\T>\nu]$ as a function of $\nu\ge0$ for the two procedures in consideration. \autoref{fig:SRP_SR_r_ADDk_vs_k_ARL_100} depicts how the sequence $\ADD_{\nu}(\T)$, indexed by $\nu$, evolves as $\nu$ runs from $0$ to $20$ for the SR--$r$ procedure (with $R_0^r=r^{*}\approx2$) and for the SRP procedure. It can be seen that $\ADD_0(\mathcal{S}_A^r)\approx\ADD_{\tinyinfty}(\mathcal{S}_A^r)$, as we planned. More importantly, note that the SR--$r$ procedure is uniformly (i.e., for all $\nu\ge0$) better than the SRP rule, while the difference is small. Starting the SR--$r$ procedure from the point that equates the average detection delays at zero and at infinity is practically more convenient, as it does not require one to know the lower bound (not to mention the quasi-stationary distribution). As this example illustrates, it may also be sufficient to outperform the SRP procedure (though for this example the gain is practically negligible).

We now turn to the accuracy of the asymptotic approximations for the average detection delays~\eqref{approxADD}. According to these approximations for both procedures the worst ADD is about $2.9$ (note that both procedures have the same threshold). However, the actual ADD-s are $3.54$ for the SRP procedure and $3.52$ for the SR--$r$ procedure. Hence, the approximations are not too accurate, which is because the ARL to false alarm is only 100.

\begin{figure}
   \centering
   \includegraphics[width=0.95\textwidth]{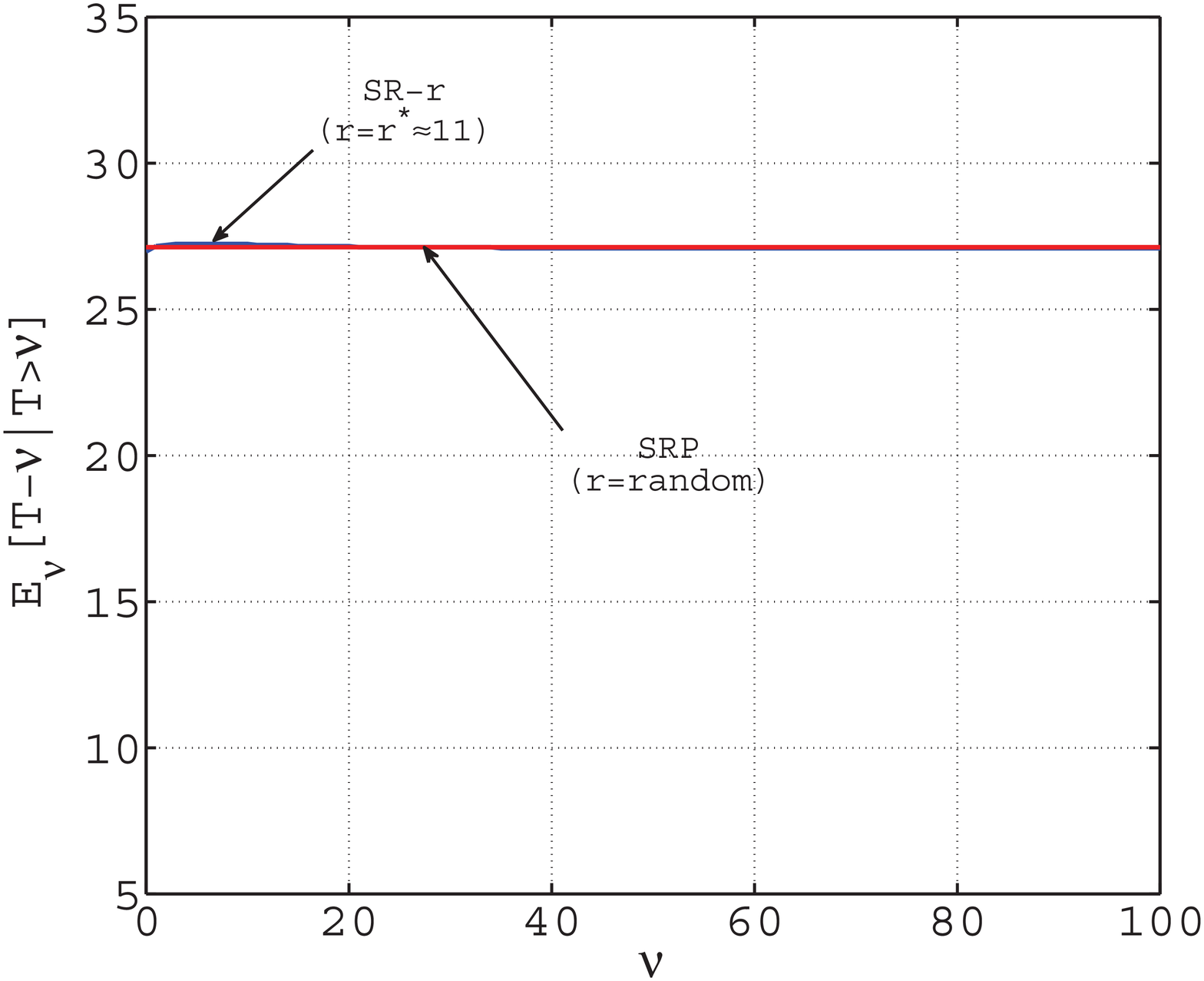}
   % figure caption to be below the figure
   \caption{Conditional average detection delay
$\EV_{\nu}[\mathcal{S}_A^r-\nu|\mathcal{S}_A^r>\nu]$ vs. change-point
$\nu$ for the SRP procedure and for the SR--$r$ procedure with
$r=r^{*}\approx11$ for the
$\mathsf{beta}(\delta,\delta+1)$-to-$\mathsf{beta}(\delta+1,\delta)$
model with $\delta=5$. The ARL to false alarm
$\EV_{\tinyinfty}[\T]=\gamma$ is approximately $5\times10^3$ for each
procedure.}
   \label{fig:SRP_SR_r_ADDk_vs_k_ARL_5000}
\end{figure}

Consider now setting $\delta$ to $5$. Since $I=1/\delta$ this is a less contrast change than $\delta=1$. Consequently, the ADD-s should be higher, which can be used to better illustrate the accuracy of their respective approximations. For $\delta=5$, we have $I=0.2$, $C_{\tinyinfty}\approx3.19$, $\zeta\approx0.685$, $\varkappa\approx0.435$, and $r^*\approx11$. Let $\gamma=5\times10^3$. To have this level of the ARL to false alarm, the threshold for the SR--$r$ procedure should be set to $3452$ (the actual ARL to false alarm for this threshold is $4999.3$), and for the SRP procedure -- to $3462$ (the actual ARL to false alarm for this threshold is $5000.1$, and $\mu_Q\approx26.1$). Again, both approximations $\EV_{\tinyinfty}[\mathcal{S}_A^r]\approx A/\zeta-r$ and $\EV_{\tinyinfty}[\mathcal{S}_A^Q]\approx A/\zeta-\mu_Q$ are highly accurate.
We now look at the delays. \autoref{fig:SRP_SR_r_ADDk_vs_k_ARL_5000} shows the average delay to detection $\ADD_{\nu}(\T)$ versus the changepoint $\nu$  for the SR--$r$ procedure with $R_0^r=r^{*}\approx11$ and for the SRP procedure. It can be seen that again $\ADD_0(\mathcal{S}_A^r)\approx\ADD_{\tinyinfty}(\mathcal{S}_A^r)$. Furthermore, the SR--$r$ procedure is almost an equalizer: there is a tiny mound raising above the SRP's flat line, though the mound is comparable in magnitude to the numerical error, and therefore, can be disregarded from a practical point of view. Both procedures are equally efficient, but since the SR--$r$ procedure is easier to initialize it is preferable for practical purposes.

In terms of the accuracy the actual $\ADD_0(\mathcal{S}_A^r)$ is $27$, while that predicted by the approximation is $27$. The actual value of $\ADD_{\tinyinfty}(\mathcal{S}_A^r)$ is $27.1$ versus the approximated value $27$ (which is the same as the value predicted for $\ADD_0(\mathcal{S}_A^r)$, because $C(r^*)=C_{\tinyinfty}$). Lastly, for the SRP procedure the actual average delay is $27.1$, while the predicted using the asymptotic approximation value is $27$. As we can see, the approximations for the ADD-s are now accurate. The reason is that the ARL to false alarm is relatively high.

To draw a line under this example, the main conclusion is that the SR--$r$ procedure is almost equalizer, and its performance is almost indistinguishable from that of the SRP procedure. However, it is easier to implement in practice, which is contrary to the SRP procedure. Hence, we recommend the SR--$r$ procedure for practical purposes.

%-------------------------------------------------------------------------------------------------%
\subsection{Example 2: An exponential scenario}

Suppose the sequence $\{X_n\}_{n\ge1}$ is comprised by the exponentially distributed random variables that undergo a shift in the mean from $1$ to $1+\theta$, where $\theta>0$. Formally, the pre- and post-change densities in this case are
\begin{align*}
f(x)
&=
\exp\left\{-x\right\}\indicator{x\ge0}
\;\;\text{and}\;\;
g(x)=\frac{1}{1+\theta}
\exp\left\{-\frac{x}{1+\theta}\right\}\indicator{x\ge0},
\end{align*}
respectively. We refer to this model as the $\mathcal{E}(1)$-to-$\mathcal{E}(1+\theta)$ model.

This model was considered by~\cite{Tartakovskyetal-IWSM09} for $\theta=0.1$, which corresponds to a small, not easily detectable change. Using the numerical framework of~\cite{Moustakidesetal-SS11}, also presented in~\autoref{s:num-perf-eval}, they carried out a performance analysis of CUSUM, the SRP procedure and the SR--$r$ procedure comparing each against the other. They also computed the lower bound. We present an excerpt of results for the SRP and SR--$r$ procedures along with the lower bound. The accuracy is within $0.5\%$.

\begin{figure}[p]
    \centering
    \includegraphics[width=0.85\textwidth]{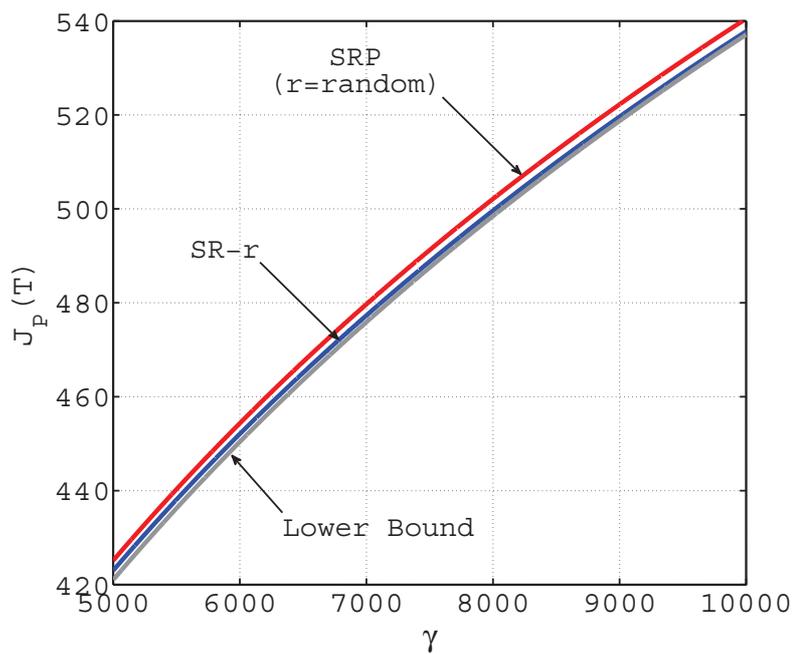}
    % figure caption to be below the figure
    \caption{The lower bound $\mathcal{J}_{\mathrm{B}}(\T)$ and Pollak's $\mathcal{J}_{\mathrm{P}}(\T)$ for the SRP and SR--$r$ procedures for the $\mathcal{E}(1)$-to-$\mathcal{E}(1+\theta)$ model with $\theta=0.1$. The ARL to false alarm is between $5\times10^3$ and $10^4$.}
    \label{fig:exp-exp-scenario:SADD-vs-ARL}
\end{figure}

\autoref{fig:exp-exp-scenario:SADD-vs-ARL} shows operating characteristics in terms of Pollak's supremum conditional average detection delay $\mathcal{J}_{\mathrm{P}}(\T)=\sup_{\nu}\EV_{\nu}[\T-\nu|\T>\nu]$ as a function of the ARL to false alarm $\EV_{\tinyinfty}[\T]=\gamma$, plus the lower bound $\mathcal{J}_{\mathrm{B}}(\T)$. It can be seen that the best performance is delivered by the SR--$r$ procedure. This is expected since by design the SR--$r$ rule is the closest to the lower bound $\mathcal{J}_{\mathrm{B}}(\T)$. This suggests that the (unknown) optimal procedure can offer only a practically insignificant improvement over the SR--$r$ procedure.

\begin{figure}[p]
    \centering
    \includegraphics[width=0.85\textwidth]{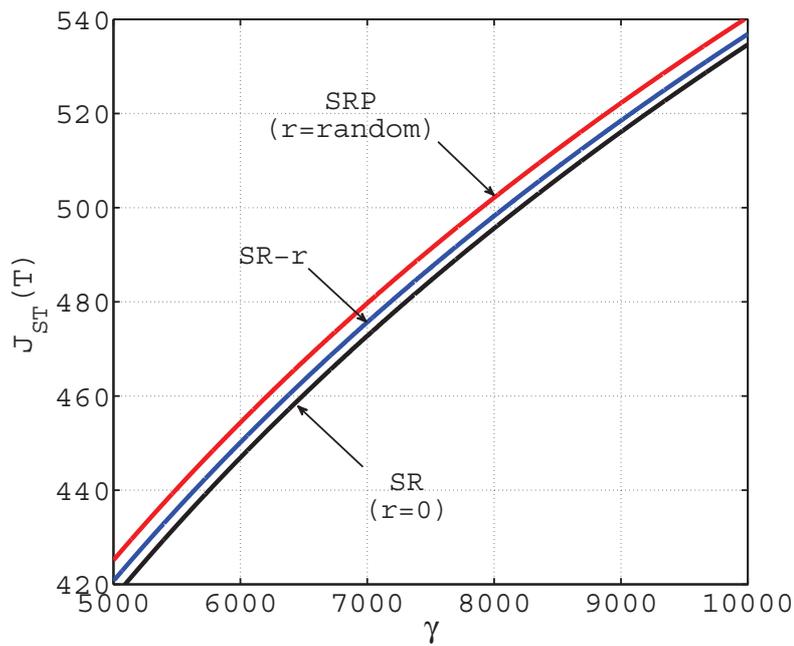}
    % figure caption to be below the figure
    \caption{The stationary average detection delay $\mathcal{J}_{\mathrm{ST}}(\T)$ for the SRP and SR--$r$ procedures for the $\mathcal{E}(1)$-to-$\mathcal{E}(1+\theta)$ model with $\theta=0.1$. The ARL to false alarm is between $5\times10^3$ and $10^4$.}
    \label{fig:exp-exp-scenario:STADD-vs-ARL}
\end{figure}

Next,~\autoref{fig:exp-exp-scenario:STADD-vs-ARL} shows the behavior of the stationary average detection delay $\mathcal{J}_{\mathrm{ST}}(\T)$ against the ARL to false alarm. Since the SR procedure is exactly optimal with respect to $\mathcal{J}_{\mathrm{ST}}(\T)$ its performance is the best among the three procedures, but the difference is relatively small. Note also that for the SRP procedure $\mathcal{J}_{\mathrm{P}}(\mathcal{S}_A^Q)$ is the same as $\mathcal{J}_{\mathrm{ST}}(\mathcal{S}_A^Q)$, since the SRP procedure is an equalizer.

%\vspace{-0.5cm}

%-------------------------------------------------------------------------------------------------%
\bibliographystyle{spbasic}      % basic style, author-year citations
\bibliography{main}

\end{document}